\documentclass[12pt, draft]{article}
\usepackage{amsmath}
\usepackage{amssymb}
\usepackage{latexsym}
\allowdisplaybreaks[4]
\newtheorem{Thm}{Theorem}[section]
\newtheorem{Prop}[Thm]{Proposition}
\newtheorem{Lemma}[Thm]{Lemma}

\newcommand{\N}{{\rm {\bf N}}}
\newcommand{\R}{{\rm {\bf R}}}
\newcommand{\C}{{\rm {\bf C}}}
\def\Bbb R{{\rm \bf R}}
\renewcommand{\v}[1]{\vert #1 \vert}

\title{A study on finding a buried obstacle in a layered medium 
having the influence of the total reflection phenomena 
via the time domain enclosure method}
\author{
Masaru IKEHATA\thanks{
Laboratory of Mathematics,
Graduate School of Engineering,
Hiroshima University, Higashihiroshima 739-8527, JAPAN;
ikehata@hiroshima-u.ac.jp
}
\and 
Mishio KAWASHITA\thanks{Department of Mathematics,
Graduate School of Sciences,
Hiroshima University, Higashihiroshima 739-8526, JAPAN;
kawasita@hiroshima-u.ac.jp}
\and
Wakako KAWASHITA\thanks{
Laboratory of Mathematics,
Graduate School of Engineering,
Hiroshima University, Higashihiroshima 739-8527, JAPAN;
wakawa@hiroshima-u.ac.jp}
}
\date{}
\begin{document}
\maketitle
\begin{abstract}
An inverse obstacle problem for the wave governed by the wave equation in a two layered medium 
is considered under the framework of the time domain enclosure method. 
The wave is generated by an initial data supported on a closed ball in the upper half-space, 
and observed on the same ball over a finite time interval. 
The unknown obstacle is penetrable and embedded in the lower half-space. 
It is assumed that the propagation speed of the wave in the upper half-space is greater than 
that of the wave in the lower half-space, which is excluded in the previous study: 
Ikehata and Kawashita (2018) to appear, Inverse Problems and Imaging. 
In the present case, when the reflected waves from the obstacle enter the upper layer, 
the total reflection phenomena occur, which give singularities to the integral representation of 
the fundamental solution for the reduced transmission problem in the background medium. 
This fact makes the problem more complicated. 
However, it is shown that these waves do not have any influence on the leading profile 
of the indicator function of the time domain enclosure method.
\end{abstract}
\vskip 1truepc
\setcounter{section}{0}
\section{Introduction and the statement of the result}\label{Introduction}
\vskip0pt\noindent

Continued on \cite{transmission No1}, 
we pursue further study on an inverse obstacle problem for the wave governed 
by a scalar wave equation in a two layered medium 
under the framework of the {\it time domain enclosure method} 
\cite{I4, IE, IEO2, IWALL, IR}.  
It is a mathematical formulation of a typical and important inverse obstacle problem and 
the solution may give us a hint to treat other inverse obstacle problems using 
electromagnetic waves, e.g., those coming from application to subsurface radar, 
ground probing radar \cite{DGS} and through-wall imaging \cite{BA}.

\par

In \cite{transmission No1} it is assumed that the unknown obstacle is penetrable and 
embedded in the lower half-space, and that the propagation speed of the wave 
in the upper half-space is {\it less} than that of the wave in the lower half-space.
The wave is generated by an initial data supported on an open ball
in the upper half-space and observed on the same ball
over a finite time interval.
It is shown that one can extract 
the {\it optical distance} from the ball to the obstacle and its qualitative property 
from the leading profile of the {\it indicator function}, which 
can be computed by using the wave observed over a finite time interval.

\par

When the propagation speed of 
the wave in the upper half-space is {\it greater} that of the wave in the lower half-space, 
the {\it total reflection} phenomena of the reflected wave 
by the obstacle may occur and complicate the problem more.  
The purpose of this article is to show that 
the leading profile of the indicator function is the same as the case 
treated in \cite{transmission No1}.  

Let $0 < T < \infty$.
Given $f \in L^2(\Bbb R^3)$ let $u=u(x,t)$ be the solution 
of the following initial value problem: 
\begin{equation}
\left\{
\begin{array}{ll}
\displaystyle
(\partial_t^2-\nabla\cdot\gamma\nabla) u = 0 & \text{in}\, (0, T)\times\Bbb R^3,
\\
\displaystyle
u(0, x) = 0, \quad
\partial_tu(0, x) = f(x) & \text{on}\,\Bbb R^3,
\end{array}
\right.
\label{transmission equation}
\end{equation}
where $\gamma=\gamma(x)=(\gamma_{ij}(x))$ satisfies
\par
$\bullet$  for each $i,j=1,2,3$ $\gamma_{ij}(x)=\gamma_{ji}(x)\in L^{\infty}(\Bbb R^3)$;
\par
$\bullet$  there exists a positive constant $C$ such that 
$\gamma(x)\xi\cdot\xi\ge C\vert\xi\vert^2$ 
for all $\xi\in\Bbb R^3$ and a.e. $x\in\Bbb R^3$.  
\par

As given in \cite{DL} (see e.g. Theorem 1 on p. 558 of \cite{DL}), 
for $f \in L^2(\Bbb R^3)$, there exists a unique 
$u \in L^2(0, T; H^1(\Bbb R^3))$ with 
$\partial_tu \in L^2(0, T; H^1(\Bbb R^3))$, 
$\partial_t^2u \in L^2(0, T; (H^1(\Bbb R^3))')$, 
such that for all $\phi \in H^1(\R^3)$, $u$ satisfies 
$$
\langle \partial_t^2u(t, \cdot), \phi \rangle
+ \int_{\Bbb R^3}\gamma(x)\nabla_xu(t, x)\cdot\nabla_x\phi(x)dx = 0 \quad
\text{a.e. }  t \in (0, T)
$$
and $u(0, x) = 0$ and $\partial_tu(0, x) = f(x)$.
This function $u$ is called the (weak) solution of $u$ of (\ref{transmission equation}).
\par

As a background medium we choose the whole space $\Bbb R^3$ and divide the space
into two homogeneous and isotropic media: 
$$\displaystyle
\Bbb R^3=\overline{\Bbb R^3_{+}}\cup\overline{\Bbb R^3_{-}},
$$
where $\Bbb R^3_{\pm} = \{ x = (x_1,x_2,x_3) \in \Bbb R^3 \,\vert\, \pm x_3 > 0 \}$.
The propagation speed of the wave in $\Bbb R^3_{\pm}$
is given by $\sqrt{\gamma_\pm}$, where $\gamma_\pm > 0$ are constants. 
We call $\Bbb R^3_+$ 
(resp. $\Bbb R^3_{-}$) the upper (resp. lower) side of the flat transmission boundary 
$\partial\Bbb R^3_{\pm}$.

\par

Now we specify the form of $\gamma$ in (\ref{transmission equation}).
Let $D$ be a bounded open set with $C^2$ boundaries satisfying 
$\overline{D} \subset \Bbb R^3_-$. 
We assume that $\gamma$ takes the form 
$$
\gamma(x) = \left\{
\begin{array}{ll}
\displaystyle \gamma_0(x)I_3, & \quad\text{if $x \in \Bbb R^3\setminus D$,}
\\
\displaystyle \gamma_0(x)I_3+h(x), & \quad\text{if $x \in D$, }
\end{array}
\right.
$$
where $\gamma_0(x) = \gamma_\pm$ for ${\pm}x_3 > 0$ and $h(x) = (h_{ij}(x)) \in L^\infty(D)$.
\vskip1pc
\hskip5mm
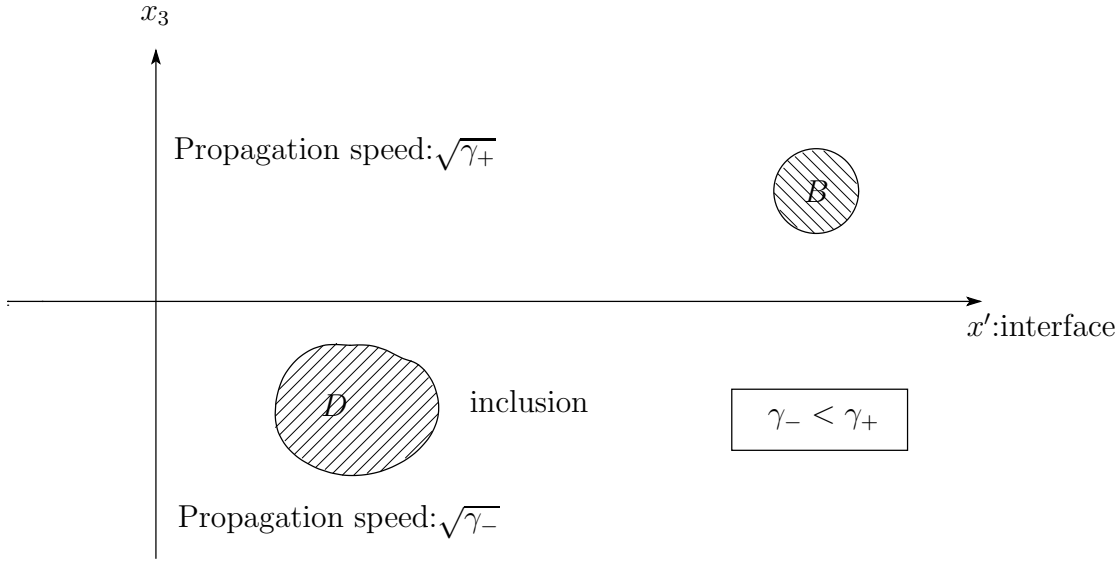
\begin{figure}[h]
\unitlength 0.1in
\unitlength 0.1in
\begin{picture}( 50.9000, 29.2500)(  3.4000,-35.5500)
%
\special{pn 8}%
\special{pa 2132 2440}%
\special{pa 2100 2438}%
\special{pa 2068 2436}%
\special{pa 2036 2436}%
\special{pa 2004 2438}%
\special{pa 1974 2446}%
\special{pa 1942 2456}%
\special{pa 1914 2470}%
\special{pa 1886 2488}%
\special{pa 1860 2508}%
\special{pa 1836 2532}%
\special{pa 1816 2556}%
\special{pa 1796 2584}%
\special{pa 1780 2612}%
\special{pa 1768 2642}%
\special{pa 1758 2672}%
\special{pa 1750 2704}%
\special{pa 1744 2736}%
\special{pa 1742 2768}%
\special{pa 1740 2800}%
\special{pa 1742 2834}%
\special{pa 1748 2866}%
\special{pa 1756 2896}%
\special{pa 1770 2924}%
\special{pa 1786 2952}%
\special{pa 1806 2978}%
\special{pa 1828 3000}%
\special{pa 1852 3022}%
\special{pa 1880 3042}%
\special{pa 1908 3058}%
\special{pa 1938 3074}%
\special{pa 1968 3088}%
\special{pa 2000 3098}%
\special{pa 2030 3108}%
\special{pa 2062 3114}%
\special{pa 2094 3120}%
\special{pa 2126 3122}%
\special{pa 2160 3122}%
\special{pa 2192 3120}%
\special{pa 2224 3116}%
\special{pa 2256 3110}%
\special{pa 2288 3102}%
\special{pa 2322 3090}%
\special{pa 2352 3078}%
\special{pa 2384 3062}%
\special{pa 2414 3046}%
\special{pa 2444 3026}%
\special{pa 2470 3006}%
\special{pa 2496 2982}%
\special{pa 2518 2958}%
\special{pa 2540 2932}%
\special{pa 2558 2904}%
\special{pa 2572 2876}%
\special{pa 2584 2846}%
\special{pa 2592 2816}%
\special{pa 2596 2784}%
\special{pa 2596 2750}%
\special{pa 2592 2718}%
\special{pa 2584 2686}%
\special{pa 2574 2654}%
\special{pa 2560 2624}%
\special{pa 2542 2596}%
\special{pa 2522 2570}%
\special{pa 2498 2548}%
\special{pa 2472 2530}%
\special{pa 2442 2520}%
\special{pa 2410 2510}%
\special{pa 2382 2496}%
\special{pa 2354 2482}%
\special{pa 2324 2466}%
\special{pa 2294 2456}%
\special{pa 2264 2446}%
\special{pa 2232 2440}%
\special{pa 2200 2438}%
\special{pa 2168 2438}%
\special{pa 2136 2440}%
\special{pa 2132 2440}%
\special{sp}%
%
\special{pn 4}%
\special{pa 2282 2450}%
\special{pa 1792 2940}%
\special{fp}%
\special{pa 2322 2470}%
\special{pa 1812 2980}%
\special{fp}%
\special{pa 2362 2490}%
\special{pa 1842 3010}%
\special{fp}%
\special{pa 2402 2510}%
\special{pa 1882 3030}%
\special{fp}%
\special{pa 2452 2520}%
\special{pa 1912 3060}%
\special{fp}%
\special{pa 2492 2540}%
\special{pa 1952 3080}%
\special{fp}%
\special{pa 2522 2570}%
\special{pa 2002 3090}%
\special{fp}%
\special{pa 2542 2610}%
\special{pa 2042 3110}%
\special{fp}%
\special{pa 2562 2650}%
\special{pa 2092 3120}%
\special{fp}%
\special{pa 2582 2690}%
\special{pa 2152 3120}%
\special{fp}%
\special{pa 2592 2740}%
\special{pa 2222 3110}%
\special{fp}%
\special{pa 2592 2800}%
\special{pa 2292 3100}%
\special{fp}%
\special{pa 2562 2890}%
\special{pa 2422 3030}%
\special{fp}%
\special{pa 2232 2440}%
\special{pa 1762 2910}%
\special{fp}%
\special{pa 2172 2440}%
\special{pa 1742 2870}%
\special{fp}%
\special{pa 2112 2440}%
\special{pa 1742 2810}%
\special{fp}%
\special{pa 2062 2430}%
\special{pa 1742 2750}%
\special{fp}%
\special{pa 1992 2440}%
\special{pa 1762 2670}%
\special{fp}%
%
\special{pn 8}%
\special{pa 1116 3556}%
\special{pa 1116 896}%
\special{fp}%
\special{sh 1}%
\special{pa 1116 896}%
\special{pa 1096 962}%
\special{pa 1116 948}%
\special{pa 1136 962}%
\special{pa 1116 896}%
\special{fp}%
\put(20.5200,-27.5000){\makebox(0,0){$D$}}%
\put(45.7200,-16.3200){\makebox(0,0){$B$}}%
\put(57.5000,-23.3000){\makebox(0,0){$x'$:interface}}%
\put(11.1500,-7.1500){\makebox(0,0){$x_3$}}%
\put(12.1000,-13.5000){\makebox(0,0)[lt]{Propagation speed:$\sqrt{\gamma_+}$}}%
\put(12.3000,-34.5000){\makebox(0,0)[lb]{Propagation speed:$\sqrt{\gamma_-}$}}%
%
\special{pn 8}%
\special{pa 520 2210}%
\special{pa 520 2210}%
\special{fp}%
%
\special{pn 8}%
\special{pa 340 2230}%
\special{pa 340 2230}%
\special{fp}%
\special{pa 340 2210}%
\special{pa 340 2210}%
\special{fp}%
\special{pa 340 2210}%
\special{pa 5430 2210}%
\special{fp}%
\special{sh 1}%
\special{pa 5430 2210}%
\special{pa 5364 2190}%
\special{pa 5378 2210}%
\special{pa 5364 2230}%
\special{pa 5430 2210}%
\special{fp}%
\put(30.7000,-27.3000){\makebox(0,0){inclusion}}%
\put(46.1000,-28.2000){\makebox(0,0){$\gamma_- < \gamma_+$}}%
%
\special{pn 8}%
\special{pa 4130 2670}%
\special{pa 5050 2670}%
\special{pa 5050 2990}%
\special{pa 4130 2990}%
\special{pa 4130 2670}%
\special{fp}%
%
\special{pn 8}%
\special{ar 4572 1632 222 222  0.0000000 6.2831853}%
%
\special{pn 4}%
\special{pa 4612 1842}%
\special{pa 4362 1592}%
\special{fp}%
\special{pa 4662 1832}%
\special{pa 4372 1542}%
\special{fp}%
\special{pa 4702 1812}%
\special{pa 4402 1512}%
\special{fp}%
\special{pa 4732 1782}%
\special{pa 4422 1472}%
\special{fp}%
\special{pa 4752 1742}%
\special{pa 4462 1452}%
\special{fp}%
\special{pa 4772 1702}%
\special{pa 4502 1432}%
\special{fp}%
\special{pa 4792 1662}%
\special{pa 4552 1422}%
\special{fp}%
\special{pa 4782 1592}%
\special{pa 4612 1422}%
\special{fp}%
\special{pa 4562 1852}%
\special{pa 4352 1642}%
\special{fp}%
\special{pa 4472 1822}%
\special{pa 4382 1732}%
\special{fp}%
\end{picture}%
\caption{Setting of the problem}
\label{Gainennzu}
\end{figure}
\vskip1pc

Note that $D$ is a mathematical model of a penetrable obstacle (inclusion) embedded 
in the lower half-space.
We introduce a {\it jump condition} of $\gamma(x)$ from $\gamma_0(x)I$ on $D$:
\par\noindent
$$
\text{(A)}_\pm
\qquad
\left\{
\begin{array}{ll}
\text{there exists a positive constant $C'$ such that }
\\
\text{${\pm}h(x)\xi\cdot\xi \geq C'\vert\xi\vert^2$ ($\xi \in \Bbb R^3$ and a.e. $x\in D$).}
\end{array}
\right.
$$

\par

We consider the following problem:
\par\noindent 
{\bf\noindent Problem.}
Fix a large $T$ (to be determined later).
Assume that $\gamma_0$ is {\it known}, $\gamma_{+}\not=\gamma_{-}$ and that both $D$ and $h$ 
are {\it unknown}. Let $B$ be an open ball with $\overline{B} \subset \Bbb R_{+}^3$.
Fix some $f\in L^2(\Bbb R^3)$ satisfying $\text{supp}\,f\subset\overline B$ and 
${\rm ess.inf}_{x \in B}f(x) > 0$ (or $-{\rm ess.inf}_{x \in B}f(x) > 0$).  
Extract information about the location and shape of $D$ from the measured data $u$ on $B$
over the time interval $(0,\,T)$, where $u$ is the weak solution of (\ref{transmission equation}) 
for the above $f$. 

\par

Note that the problem asks us to extract some information about the unknown obstacle from a 
{\it single} observed wave over a {\it finite time} interval.  The place where the wave is 
observed is the same as the generating place of the wave.  This is a near field version of 
the inverse backscattering problem in the {\it time domain}
and different from the studies in \cite{DEKPS, LZ, LLLL} where 
the {\it time harmonic reduced case} in a two layered medium have been treated.

\par

Note also that the case where $\gamma_{+}=\gamma_{-}$ has been considered in \cite{IE} and applying
an idea in \cite{IEO2} to this case, one can extract the distance 
of the ball $B$ to the obstacle $D$, that is 
${\rm dist}(D, B) = \inf_{x \in D, y \in B}\vert x - y \vert$.
Moreover, a similar inverse obstacle problem for the wave governed by the equation 
$(\alpha(x)\partial_t^2-\Delta)u=0$ has been considered in \cite{IWALL} for a general 
inhomogeneous background medium. 
In this case lower and upper estimates of ${\rm dist}(D, B)$ are given.

\par

To describe a solution to the present problem we recall the definition of 
the {\it optical distance} between the ball $B$ and obstacle $D$ given by
\begin{equation*}
\displaystyle
l(D,B)=\inf_{x\in D,\,y\in B}l(x, y),
\end{equation*}
where
\begin{align}
l(x, y) &= \inf_{z' \in \Bbb R^2}l_{x, y}(z'),
\label{definition of l(x, y)}
\\
l_{x, y}(z') &= \frac{1}{\sqrt{\gamma_-}}\v{\tilde{z}'-x}
+\frac{1}{\sqrt{\gamma_+}}\v{\tilde{z}'-y}
\quad(\tilde{z}' = (z_1, z_2, 0), z' = (z_1, z_2)).
\label{the path length for Snell's law}
\end{align}
As is in Lemma 4.1 of \cite{transmission No1}, for arbitrary $x$ and $y \in \R^3$ 
with $x_3 < 0$ and $y_3 > 0$, there exists the unique point $z'(x, y) \in \R^2$ 
satisfying $l(x, y) = l_{x, y}(z'(x, y))$, and the point $z'(x, y)$ is on the
line segment $x'y'$ and $C^\infty$ for $x$ and $y \in \R^3$ with 
$x_3 < 0$ and $y_3 > 0$. 
\par

Recall the indicator function given in \cite{transmission No1}:
$$
I_f(\tau,T) = \int_{\Bbb R^3}f(x)(w(x, \tau)- v(x, \tau))dx,
$$
where 
$$
w(x, \tau) = \int_0^Te^{-{\tau}t}u(t, x)dt \qquad(x \in \R^3)
$$
and $v \in H^1(\Bbb R^3)$ is the weak solution of 
\begin{equation*}
\begin{array}{ll}
\displaystyle
(\nabla\cdot\gamma_0\nabla-\tau^2)v(x, \tau) + f(x) = 0 & 
\displaystyle
\text{in}\,\Bbb R^3.
\end{array}
\end{equation*}

\par

Our main result is the following theorem:

\begin{Thm}\label{the goal of this article} Assume that $\gamma_{+} > \gamma_{-}$, then we have;
$$
\displaystyle
\lim_{\tau\rightarrow\infty}e^{\tau T}I_{f}(\tau, T)
=
\left\{
\begin{array}{ll}
0, & \text{if $T<2l(D,B)$,}
\\
\mp\infty,
&
\text{if $T>2l(D,B)$ and $\gamma$ satisfies $\text{(A)}_\pm$.}
\end{array}
\right.
$$
Moreover, if $\gamma$ satisfies $\text{(A)}_\pm$, then for all $T>2l(D,B)$
\begin{equation}
\lim_{\tau\longrightarrow\infty}\frac{1}{\tau}
\log\left\vert I_{f}(\tau,T) \right\vert 
= -2l(D,B).
\label{the shortest time appeared!!}
\end{equation}
\end{Thm}
Note that we have obtained the same result as \cite{transmission No1}.  
Our studies have completely covered the case $\gamma_{+}\not=\gamma_{-}$.
In the case $\gamma_{+}<\gamma_{-}$, the incident waves from the lower half-space 
do not cause the total reflection. 
On the other hand, in the present 
case $\gamma_+ > \gamma_{-}$ those waves cause the total reflection, 
which makes the problems more complicated than that 
of \cite{transmission No1}. Theorem \ref{the goal of this article} shows that
the total reflection phenomena do not have any influence on the leading profile
of the indicator function $I_{f}(\tau,T)$ as $\tau \to \infty$.

\par

Theorem \ref{the goal of this article} says that we need to take $T > 2l(D, B)$ at least
if we wish to know information of $D$ from the indicator function. 
We think this restriction is optimal and consistent with wave phenomena, 
since we should wait to the signals going and coming back to the points taking measurements.
We can also know whether the propagation speed 
of wave in the inclusion is greater or less than the speed of wave in the background medium
by checking the asymptotic behavior of $e^{\tau T}I_f(\tau,T)$ as $\tau\rightarrow\infty$.
From the formula (\ref{the shortest time appeared!!}), we can compute the value $l(D, B)$. 
Moreover, as pointed out in \cite{transmission No1}
we have
$$\displaystyle
D\subset E(D;B,\gamma_+,\gamma_{-}),
$$
where
$$\displaystyle
E(D;B,\gamma_+,\gamma_{-})
=\left\{x\in\Bbb R_{-}^3\,\mid l(x,p)>l(D,B)+\frac{\eta}{\sqrt{\gamma_+}}\right\}
$$
and $p$ and $\eta$ are the center point and radius of $B$, respectively.  Note that
the set $E(D;B,\gamma_+,\gamma_{-})$ can be determined
by the computed value of $l(D,B)$, $\eta$ and $\sqrt{\gamma_{+}}$.
This means that the one shot yields one information about the geometry of $D$.

\par

The proof of Theorem \ref{the goal of this article} proceeds along the same lines as the case
$\gamma_+ < \gamma_-$ in \cite{transmission No1}. 
The indicator function has the well known estimates 
below (see Lemma 1.2 in \cite{transmission No1}).
\begin{Lemma}\label{relations between the indicator function and v}
We have, as $\tau\longrightarrow\infty$
\begin{align}
I_f(\tau,T) &\geq 
\int_{\Bbb R^3}(\gamma_0I_3-\gamma)\nabla v\cdot\nabla v dx+O(\tau^{-1}e^{-\tau T})
\label{for (A)_-}
\intertext{and}
I_f(\tau,T) &\leq \int_{\Bbb R^3}\gamma_0(\gamma_0I_3-\gamma)\gamma^{-1/2}
\nabla v\cdot\gamma^{-1/2}\nabla vdx+O(\tau^{-1}e^{-\tau T}).
\label{for (A)_+}
\end{align}
\end{Lemma}

From (\ref{for (A)_-}) and $\text{(A)}_-$ (resp. (\ref{for (A)_+}) and $\text{(A)}_+$), 
we  see that Theorem \ref{the goal of this article} immediately follows from the following 
estimates for $v$. 
\begin{Thm}\label{estimate of nabla_xv on D} 
Assume that $\partial D$ is $C^1$ and that $\gamma_{+} > \gamma_{-}$.
Then, there exist positive numbers $C$ and $\tau_0$ such that, for all $\tau\ge\tau_0$
we have
\begin{equation*}
C^{-1}\tau^{-4}e^{-2\tau l(D,B)} \leq \int_{D}\vert\nabla v(x)\vert^2\,dx 
\leq C\tau^2e^{-2\tau l(D,B)}.
\end{equation*}
\end{Thm}
Note that Theorem \ref{estimate of nabla_xv on D} in which
the assumption $\gamma_{+}>\gamma_{-}$ is replaced with $\gamma_{+}<\gamma_{-}$ 
has been established in \cite{transmission No1}. 

\par 

Thus everything is reduced to showing the validity of Theorem \ref{estimate of nabla_xv on D}.
For the purpose we make use of the expression
$$
v(x) = \int_{B}\Phi_\tau(x, y)f(y)dy,
$$
where $\Phi_\tau(x, y)$ is governed by 
$$
\nabla_x\cdot(\gamma_0(x)\nabla_x \Phi_\tau(x, y)) -\tau^2\Phi_\tau(x, y)
+\delta(x-y) = 0
\qquad \text{in } \Bbb R^3.
$$
Since it follows that
\begin{align*}
\int_{D}\v{\nabla_xv(x)}^2dx
= \int_{B}dy\int_{B}d{\xi}f(y){f(\xi)}
\int_{D}\nabla_x\Phi_\tau(x, y)\cdot{\nabla_x\Phi_\tau(x, \xi)}dx, 
\end{align*}
Theorem \ref{estimate of nabla_xv on D} is given by investigating 
an asymptotic behavior of $\nabla_x\Phi_\tau(x, y)$
as $\tau \to \infty $ for $x = (x', x_3)$ with $x_3 < 0$, $x' \in \Bbb R^2$
and $y \in B$.

In section \ref{Asymptotics and estimates of the refracted part}, 
a complex integral representation of 
the fundamental solution $\Phi_\tau(x, y)$ is recalled, which is given in \cite{transmission No1}. 
As in (\ref{expression of the FS for x_3 < 0}) in 
section \ref{Asymptotics and estimates of the refracted part}, 
$\Phi_\tau(x, y)$ consists of the part 
corresponding to the incident wave and the refracted part $ E^{\gamma_-}_{\tau}(x, z') $
for $x \in \R^3_-$ with $x_3 < 0$ and $z' \in \R^2$. 
To obtain asymptotics for $\Phi_\tau(x, y)$, the steepest descent method is 
used for the integral representation of the refracted part. 
If $\gamma_+ < \gamma_-$, the integrand in the representation of the refracted part
is holomorphic near the steepest descent curve. 
Hence, we can perform asymptotic expansion of the refracted part 
and $\Phi_\tau(x, y)$.

On the other hand, if $\gamma_+ > \gamma_-$, the total reflection phenomena for incident waves
from the lower half-space occur, which correspond to the fact that the steepest descent 
curve should be across singularities of the integrand when the contour is changed. 
Because of singularities, it seems difficult to get asymptotics of the refracted part 
corresponding to the total reflection phenomena. Hence, we only obtain some estimates for 
the refracted part containing the total reflection phenomena, which is the purpose of 
section \ref{Asymptotics and estimates of the refracted part}.

\par

In section \ref{The shortest time and  asymptotics}, we show the following asymptotics of 
$\nabla_x\Phi_\tau(x, y)$ (and $\Phi_\tau(x, y)$):
\begin{Prop}\label{Asymptotics of the refracted part of the gradient of the FS}
Assume that $\gamma_+ > \gamma_-$. Then for $k = 0, 1$, we have 
\begin{align}
\nabla_x^k\Phi_\tau(x, y) = 
\frac{{e}^{-{\tau}l(x, y)}}{8\pi\gamma_+\gamma_-
\sqrt{{\rm det}H(x, y)}}
\Big(\frac{-\tau}{\sqrt{\gamma_-}}\Big)^k
\Big(\sum_{j = 0}^{N}{\tau^{-j}}\Phi_{j}^{(k)}(x, y) 
+ Q_{N, \tau}^{(k)}(x, y)\Big), 
\label{asymp of Phi}
\end{align}
where $H(x, y) = {\rm Hess}(l_{x, y})(z'(x, y))$ is the Hessian of $l_{x, y}$ given by
(\ref{the path length for Snell's law}) at $z' = z'(x, y)$, 
$\Phi_{j}^{(k)}(x, y)$ $(k = 0, 1)$ are 
$C^\infty$ in $\overline{D}\times\overline{B}$,
for any $N \in {\rm N}\cup\{0\}$, $Q_{N, \tau}^{(k)}(x, y)$ $(k = 0, 1)$
are continuous in $\overline{D}\times\overline{B}$
with a constant $C_N > 0$ satisfying
$$
\v{Q_{N, \tau}^{(0)}(x, y)} + \v{Q_{N, \tau}^{(1)}(x, y)} \leq C_N\tau^{-(N+1)} \qquad(x \in \overline{D}, y \in \overline{B}, \tau \geq 1).
$$
Moreover, $\Phi_{0}^{(k)}(x, y) $ $(k = 0, 1)$ are given by
\begin{align*}
\Phi_{0}^{(0)}(x, y)
&= \frac{E_{0}(x - \tilde{z}'(x, y))}{\v{x - \tilde{z}'(x, y)}\v{\tilde{z}'(x, y)-y}},
\intertext{and}
\Phi_{0}^{(1)}(x, y) 
&= \Phi_{0}^{(0)}(x, y)\frac{x - \tilde{z}'(x, y)}{\v{x - \tilde{z}'(x, y)}},
\end{align*}
where
\begin{align}
E_{0}(x - \tilde{z}') 
&= \frac{4\sqrt{\gamma_-}\v{x_3}\sqrt{a_0^2\v{x - \tilde{z}'}^2- \v{x' - z'}^2}}
{\v{x - \tilde{z}'}\big(\sqrt{a_0^2\v{x - \tilde{z}'}^2- \v{x' - z'}^2}+a_0^2\v{x_3}\big)}.
\label{the form of E_0}
\end{align}
\end{Prop}

Note that Proposition \ref{Asymptotics of the refracted part of the gradient of the FS} is
the same as Proposition 1 of \cite{transmission No1} except for the condition 
$\gamma_{+} > \gamma_{-}$. 
That means that the total reflection phenomena make no difference to asymptotics of 
$\Phi_\tau(x, y)$. 
We should consider the influence of the total reflection phenomena on the 
optical distance, which is discussed in 
section \ref{Asymptotics and estimates of the refracted part}.  
For the usual inner waves, 
the optical distance between $z'\in \R^2 = \partial \R_-^3$ and $x \in \R^3_-$ 
is given by $\v{x - \tilde{z}'}/\sqrt{\gamma_-}$. 
Hence, the optical distance between $x \in \R^3_-$ and $y \in \R^3_+$ is given by 
(\ref{definition of l(x, y)}) and (\ref{the path length for Snell's law}) if the total 
reflection phenomena do not occur.
In our case, we should pay attention to the fact that 
optical distance between $\tilde{z}'$ and $x$ 
corresponding to the total reflection phenomena 
is different from $\v{x - \tilde{z}'}/\sqrt{\gamma_-}$ 
(see (\ref{def of tilde{l}_{x, y}(z')}) in section \ref{The shortest time and asymptotics}).  
Hence, in this case, the time in which the waves travel from $x$ to $y$ via 
$\tilde{z}' \in \partial\R^3_+$ is also different from $l_{x, y}(z')$ given in 
(\ref{the path length for Snell's law}). 
But we can show that even in this case, the function $l(x, y)$ gives the optical distance 
between $x \in \R^3_-$ and $y \in \R^3_+$ (cf. Lemma \ref{the shortest length in the case}).
As is in \cite{transmission No1}, the fact that $l(x, y)$ gives the optical distance 
plays an important role to obtain 
Proposition \ref{Asymptotics of the refracted part of the gradient of the FS}. 
This is the reason why 
Proposition \ref{Asymptotics of the refracted part of the gradient of the FS} 
has the same conclusion as in \cite{transmission No1}.  
Once we obtain asymptotics in 
Proposition \ref{Asymptotics of the refracted part of the gradient of the FS},
Theorem \ref{estimate of nabla_xv on D} can be shown by the same argument as in 
\cite{transmission No1}. This is the outline of this article.

\par

\setcounter{equation}{0}
\section{Asymptotics and estimates of the refracted part}
\label{Asymptotics and estimates of the refracted part}
\vskip0pt\noindent

Let us recall an integral representation of the fundamental solution $\Phi_\tau(x, y)$ 
given in \cite{transmission No1}. 
A usual fundamental solution for the case of no transmission boundary (i.e. the case
of $\gamma_- = \gamma_+$) is of the form:
$$
E^{\gamma_+, 0}_\tau(x, y) = \frac{1}{4\pi\gamma_+}
\frac{e^{-\tau\v{x - y}/\sqrt{\gamma_+}}}{\v{x - y}}
\quad(x \neq y, \tau > 0),
$$ 
which coincides with that of defined by the Fourier integral 
\begin{equation*}
E^{\gamma_+, 0}_\tau(x, y) = \frac{1}{(2\pi)^3}\int_{\Bbb R^3}
{e}^{i\xi\cdot(x-y)}
\frac{1}{\gamma_+\xi^2+\tau^2}
d{\xi}
= \frac{\tau}{(2\pi)^3}\int_{\Bbb R^3}{e}^{i\tau\xi\cdot(x-y)}
\frac{1}{\gamma_+\xi^2+1}d\xi.
\end{equation*}
As in (11) of \cite{transmission No1}, we introduce
\begin{equation*}
E^{\gamma_-}_{\tau}(x, z')
= \frac{\tau}{(2\pi)^3}
\int_{\Bbb R^3}{e}^{i\tau\xi\cdot(x-\tilde{z}')}
\frac{1}{\gamma_-\xi^2+1}R(\sqrt{\gamma_-}\v{\xi'})d{\xi}
\qquad(x_3 < 0), 
\end{equation*}
where $\tilde{z}' = (z', 0)$ $(z' \in \Bbb R^2)$ is the point on 
the transmission boundary $\partial\Bbb R^3_\pm$ and
$R(\v{\xi'})$ is a function of $\v{\xi'}$ standing for 
the transmission coefficient given by 
\begin{align*}
R(\rho) 
= \frac{4\sqrt{\gamma_-}\sqrt{a_0^2 + \rho^2}\sqrt{1 + \rho^2}}
{\sqrt{a_0^2 + \rho^2}+a_0^2\sqrt{1 + \rho^2}}
\quad(\rho \geq 0)
\quad\text{with}\quad
a_0 = \sqrt{\frac{\gamma_-}{\gamma_+}}.
\end{align*}
Using $E^{\gamma_+, 0}_{\tau}(x, \tilde{z'})$ and $E^{\gamma_-}_{\tau}(x, z')$ we can
represent the fundamental solution $\Phi_\tau(x, y)$ for 
$y \in \R^3_+$ and $x \in \R^3_-$ as
\begin{align}
\Phi_\tau(x, y)
= \frac{\tau}{4\pi\gamma_+}\int_{\Bbb R^2}
E^{\gamma_-}_{\tau}(x, z')
\frac{e^{-\tau\v{\tilde{z}'- y}/\sqrt{\gamma_+}}}{\v{\tilde{z}'- y}}dz'.
\label{expression of the FS for x_3 < 0}
\end{align}
This is just (10) of \cite{transmission No1}. 
In what follows, as in \cite{transmission No1}, 
we call $E^{\gamma_-}_{\tau}(x, z')$ the refracted part (of the fundamental
solution $\Phi_\tau(x, y)$).

\par

Put $\Theta_k(x, z') = \frac{x_k - z_k}{\v{x'-z'}}$ $(k = 1, 2)$ and  
$\Theta_3(x, z') = \frac{x_3}{\v{x_3}}$. 
Note that (26)-(29) of \cite{transmission No1} imply that
the refracted part $E^{\gamma_-}_{\tau}(x, z')$ is expressed by 
\begin{align}
E^{\gamma_-}_{\tau}(x, z') &= \frac{\tau}{2(2\pi)^2\gamma_-^{3/2}}
\int_{\Bbb R}I_{\tilde\tau, 0}(x-\tilde{z}', \zeta_2)d\zeta_2, 
\label{For asymptotics of the refracted part 1-0}
\\
\partial_{x_k}E^{\gamma_-}_{\tau}(x, z')
&= \frac{\tau^2}{2(2\pi)^2\gamma_-^{2}}
\int_{\Bbb R}I_{\tilde\tau, k}(x-\tilde{z}', \zeta_2)d\zeta_2\Theta_k(x, z')
\qquad(k = 1, 2, 3),
\label{For asymptotics of the refracted part 1-1}
\end{align}
where for $x \in \Bbb R^3_-$, $z' \in \Bbb R^2$ and $k = 0, 1, 2, 3$, we put 
$\tilde\tau = \tau/\sqrt{\gamma_-}$, 
\begin{align}
I_{\tilde\tau, k}(x-\tilde{z}', \zeta_2) = \int_{\Bbb R}{e}^{-\tilde\tau\sqrt{1+\zeta_2^2}
(-i\v{x'-z'}\zeta_1+\v{x_3}\sqrt{1+\zeta_1^2})}Q_k(\zeta_1, \zeta_2)
\frac{d\zeta_1}{\sqrt{1+\zeta_1^2}}
\label{the integral in zeta_1}
\end{align}
and 
\begin{align*}
\left\{
\begin{array}{lll}
Q_0(\zeta_1, \zeta_2) = R\left(\sqrt{\zeta_1^2+\zeta_2^2+\zeta_1^2\zeta_2^2}\right), 
\quad
\tilde{Q}_0(\zeta_1, \zeta_2) = \sqrt{1+\zeta_2^2}Q_0(\zeta_1, \zeta_2), 
\\
Q_1(\zeta_1, \zeta_2) = Q_2(\zeta_1, \zeta_2) = i\zeta_1\tilde{Q}_0(\zeta_1, \zeta_2), 
\quad
Q_3(\zeta_1, \zeta_2) = -\sqrt{1+\zeta_1^2}\tilde{Q}_0(\zeta_1, \zeta_2). 
\end{array}
\right.
\end{align*}
We use the steepest descent method to the integrals $I_{\tilde\tau, k}(x-\tilde{z}', \zeta_2)$. 
Take $\theta$ satisfying \begin{equation}
\sin\theta = \frac{\v{x' - z'}}{\v{x - \tilde{z}'}}, 
\quad
\cos\theta = \frac{\v{x_3}}{\v{x - \tilde{z}'}}
\quad(0 \leq \theta \leq \pi/2),
\label{the definition of theta}
\end{equation}
and put $r = \v{x - \tilde{z}'}\sqrt{1+\zeta_2^2}$ and
\begin{align}
\lambda = 
\lambda(\zeta_1, x, z') = -i\sin\theta\zeta_1+\cos\theta\sqrt{1+\zeta_1^2}.
\label{change of variables}
\end{align}
Note that (\ref{change of variables}) is equivalent to 
$\zeta_1 = i\lambda\sin\theta\pm\sqrt{\lambda^2-1}\cos\theta$, which yields
\begin{align}
\zeta_1 = \zeta_1(\rho, x, z') = i\sqrt{1+\rho^2}\sin\theta+\rho\cos\theta
\qquad(\rho \in \Bbb R, x \in \Bbb R^3_-, z' \in \Bbb R^2)
\label{parametrization of Im lambda = 0}
\end{align}
by putting $\lambda = \sqrt{1+\rho^2}$ for $\lambda \geq 1$ 
(cf. (33) in \cite{transmission No1}).

\par

In the case of $\gamma_+ < \gamma_-$, the function 
\begin{align}
Q_0(\zeta_1, \zeta_2) 
= \frac{4\sqrt{\gamma_-}\sqrt{1+\zeta_2^2}\sqrt{1 + \zeta_1^2}P(\zeta_1, \zeta_2)}
{P(\zeta_1, \zeta_2)+a_0^2\sqrt{1 + \zeta_1^2}}, 
\label{the form of composition function of R}
\end{align}
where
$$
P(\zeta_1, \zeta_2) = \sqrt{\frac{a_0^2+\zeta_2^2}{1+\zeta_2^2}+\zeta_1^2} 
$$
is holomorphic for $\zeta_1 \in \C\setminus((-i\infty, -i]\cup[i, i\infty))$. 
Hence, we can change the contour of integrals (\ref{the integral in zeta_1})
to the curve $\Gamma_{x, z'}$ defined by (\ref{change of variables}). This implies 
\begin{align}
I_{\tilde{\tau}, k}(x-\tilde{z}', \zeta_2)
= \int_{\Gamma_{x, z'}}e^{-\tilde{\tau}r\lambda}
Q_k(\zeta_1, \zeta_2)\frac{d\zeta_1}{\sqrt{1+\zeta_1^2}}. 
\label{s of I}
\end{align}
Using this formula, we can obtain asymptotics of $\Phi_\tau(x, y)$ as $\tau \to \infty$.
On the contrary, in the case of $\gamma_+ > \gamma_-$, i.e. 
$a_0 = \sqrt{\gamma_-/\gamma_+} < 1$, the functions $P$ and $Q_0$ 
are holomorphic for 
$\zeta_1 \in \C\setminus((-i\infty, -ib_0(\zeta_2)]\cup[ib_0(\zeta_2), i\infty))$, 
where 
$$
b_0(\zeta_2) = \sqrt{\frac{a_0^2+\zeta_2^2}{1+\zeta_2^2}}.
$$
Thus, if $\sin\theta < a_0$, we can change the contour
to $\Gamma_{x, z'}$, however, if $\sin\theta > b_0(\zeta_2)$ we should make a detour to connect 
$\Gamma_{x, z'}$ and the branch point $\zeta_1 = ib_0(\zeta_2)$ of $P(\zeta_1, \zeta_2)$. 
This corresponds to the total reflection phenomena, 
which makes us additional arguments.

\par
In what follows, for $\delta$ with $0 < \delta < a_0^{-1}$ and $x \in \R^3_-$, we put 
$ 
{\mathcal U}_{\delta}(x) = \{\, z' \in \R^2 \,\vert\, \v{x' - z'} 
< a_0{\delta}\v{x - \tilde{z}'} \,\}. 
$
Note that $z' \in \overline{{\mathcal U}_{\delta}(x)}$ is equivalent to 
\begin{align}
\v{x' - z'} < \frac{a_0\delta}{\sqrt{1-a_0^2\delta^2}}\v{x_3}.
\label{estimate for x and z' in overline{{mathcal U}_{delta}(x)}}
\end{align}

\par

Since ${\mathcal U}_{\delta}(x) = \{\, z' \in \R^2 \,\vert\, 
\sin\theta < a_0\delta\,\}$, it follows that
$\inf\{ \v{ia_0 - \zeta_1} \,\vert\, \zeta_1 \in \Gamma_{x, z'} \,\} = a_0(1-\delta)$
for any $0 < \delta < 1$, $x \in \R^3_-$ and $z' \in {\mathcal U}_{\delta}(x)$. 
Thus, in this case, the argument for getting Proposition 2 in \cite{transmission No1} 
implies the following expansions of the refracted part:
\begin{Lemma}\label{Asymptotics for no total reflection waves}
Assume that $\gamma_+ > \gamma_-$. Then, 
for any $0 < \delta < 1$, the refracted part $E^{\gamma_-}_{\tau}(x, z')$ for $x \in \R^3_-$ and
$z' \in \overline{{\mathcal U}_{\delta}(x)}$ is expanded by 
\begin{align*}
E^{\gamma_-}_{\tau}(x, z')
= \frac{e^{-{\tau\v{x-\tilde{z}'}/\sqrt{\gamma_-}}}}{4\pi\gamma_-\v{x - \tilde{z}'}}
&\Big(\sum_{j = 0}^{N-1}E_{j}(x-\tilde{z}')
\Big(\frac{\sqrt{\gamma_-}}{\tau\v{x - \tilde{z}'}}\Big)^j
+\tilde{E}_{N}(x, z'; \tau)
\Big),
\intertext{and for $k = 1, 2, 3$, } 
\partial_{x_k}E^{\gamma_-}_{\tau}(x, z')
= \frac{-{\tau}e^{-{\tau\v{x-\tilde{z}'}/\sqrt{\gamma_-}}}}{4\pi\gamma_-^{3/2}\v{x - \tilde{z}'}}
&\Big(
\sum_{j = 0}^{N-1}G_{k, j}(x-\tilde{z}')\Big(\frac{\sqrt{\gamma_-}}{\tau\v{x - \tilde{z}'}}\Big)^j
+\tilde{G}_{k, N}(x, z'; \tau)\Big), 
\end{align*}
where $E_j(x - \tilde{z}')$, $G_{k, j}(x - \tilde{z}')$ 
($k = 1, 2, 3$ and $j = 0, 1, 2, \ldots)$ are $C^\infty$
functions for $x$ and $z'$ with $z' \in \overline{{\mathcal U}_{\delta}(x)}$. 
Here, the remainder terms
$\tilde{E}_{N}(x, z'; \tau)$ and $\tilde{G}_{k, N}(x, z'; \tau)$ $(k = 1, 2, 3)$ 
are estimated by 
$$
\v{\tilde{E}_{N}(x, z'; \tau)}+\sum_{k = 1}^3\v{\tilde{G}_{k, N}(x, z'; \tau)}
\leq C_{N, \delta}\Big(\frac{\sqrt{\gamma_-}}{\tau\v{x - \tilde{z}'}}\Big)^{N}
\quad(x \in \R^3_-, z' \in \overline{{\mathcal U}_{\delta}(x)})
$$
for some constant $C_{N, \delta} > 0$ depending only on $N \in \N$ and $\delta$.
In particular, we have
\begin{align*}
\left\{
\begin{array}{ll}
G_{k, 0}(x - \tilde{z}') = E_{0}(x - \tilde{z}')\displaystyle\frac{x_k - z_k}{\v{x - \tilde{z}'}}
\quad(k = 1, 2)\quad
\\[2mm]
G_{3, 0}(x - \tilde{z}') = E_{0}(x - \tilde{z}')\displaystyle\frac{x_3}{\v{x - \tilde{z}'}}, 
\end{array}
\right.
\end{align*}
where $E_{0}(x - \tilde{z}')$ is given in (\ref{the form of E_0}).
\end{Lemma}

\par

Thus, once $0 < \delta < 1$ is fixed, we can obtain uniform estimates of the refracted part for 
$x \in \R^3_-$ and $z' \in \overline{{\mathcal U}_{\delta}(x)}$. On the contrary, for 
$x \in \R^3_-$ and $z' \in \R^2\setminus{\mathcal U}_{\delta}(x)$, it seems to be hard to 
get asymptotics of the refracted part by the total reflection waves.
Fortunately, for our purpose, we have only to obtain the estimates 
for the refracted part. 
The main part of this section is to show these estimates.

\par

If $\theta$ is near $\theta_0$ and $\theta \leq \theta_0$, we have the following expansions:
\begin{Prop}\label{Asymptotics for near the total reflection angle case 1}
Assume that $\gamma_+ > \gamma_-$. Then, for any fixed $\delta$ with $0 < \delta < 1$,
the refracted part ${E}_{\tau}^{\gamma_-}(x, z')$ for $x \in \R^3_-$ and 
$z' \in \overline{{\mathcal U}_{1}(x)\setminus{\mathcal U}_{\delta}(x)}$ is expanded by 
\begin{align*}
E^{\gamma_-}_{\tau}(x, z')
&= \frac{e^{-{\tau\v{x - \tilde{z}'}/\sqrt{\gamma_-}}}}{4\pi\gamma_-\v{x - \tilde{z}'}}
\Big(E_{0}(x-\tilde{z}')+\tilde{E}_{0, 0}^{\gamma_-}(x, z'; \tau)
\Big), 
\\
\partial_{x_k}E^{\gamma_-}_{\tau}(x, z')
&= 
\frac{-{\tau}e^{-{\tau\v{x - \tilde{z}'}/\sqrt{\gamma_-}}}}{4\pi\gamma_-^{3/2}\v{x - \tilde{z}'}}
\Big(G_{k, 0}(x - \tilde{z}')+\tilde{E}_{k, 0}^{\gamma_-}(x, z'; \tau)\Big)
\quad(k = 1, 2, 3).
\end{align*}
In the above, $E_0$ and $G_{k, 0}$ are the functions given in 
Lemma \ref{Asymptotics for no total reflection waves}.
For the remainder terms $\tilde{E}_{k, 0}^{\gamma_-}(x, z'; \tau)$, for any 
$0 < \delta < 1$, there exists a constant $C_{\delta} > 0$ such that
\begin{align*}
\v{\tilde{E}_{k, 0}^{\gamma_-}(x, z'; \tau)} 
&\leq C_{\delta}\Big(\frac{\sqrt{\gamma_-}}
{\tau\v{x - \tilde{z}'}}\Big)^{1/4}
\quad(x \in \R^3_-, z \in \overline{{\mathcal U}_{1}(x)\setminus{\mathcal U}_{\delta}(x)}, 
k = 0, 1, 2, 3).
\end{align*}
\end{Prop}

For the case of $\theta > \theta_0$, we have the following estimates:
\begin{Prop}\label{Asymptotics for near the total reflection angle case 2}
Assume that $\gamma_+ > \gamma_-$. Then, there exists a constant $C > 0$ such that
the refracted part ${E}_{\tau}^{\gamma_-}(x, z')$ for $x \in \R^3_-$ and 
$z' \in \R^2\setminus{\mathcal U}_{1}(x)$ is estimated by 
\begin{align*}
\v{\nabla_x^k{E}_{\tau}^{\gamma_-}(x, z')} &\leq C\tau^{k}e^{-{\tau}T_{x, z'}(\theta_0)}
\qquad(x \in \overline{D}, z' \in \R^2\setminus{\mathcal U}_{1}(x), k = 0, 1),
\end{align*}
where for $x \in \R^3_-$ and $z' \in \R^2$, 
$T_{x, z'}(\alpha)$ is defined by
\begin{align}
T_{x, z'}(\alpha) 
= \frac{1}{\sqrt{\gamma_-}}\Big(\v{x_3}\cos\alpha + \v{z'-x'}\sin\alpha\Big).
\label{T(alpha)} 
\end{align}
\end{Prop}
Note that $T_{x, z'}(\alpha)$ is expressed by
\begin{equation}
T_{x, z'}(\alpha) = \frac{\v{\tilde{z}'-x}}{\sqrt{\gamma_-}}\cos(\theta-\alpha)
\label{expression of T(alpha)}
\end{equation}
by using $\theta$ defined by (\ref{the definition of theta}). In what follows, we only write
$T_{x, z'}(\alpha)$ by $T(\alpha)$ shortly.

The rest of this section is devoted to show 
Propositions \ref{Asymptotics for near the total reflection angle case 1} and 
\ref{Asymptotics for near the total reflection angle case 2}.

\noindent
Proof of Proposition \ref{Asymptotics for near the total reflection angle case 1}. 
When $z' \in \overline{{\mathcal U}_{1}(x)\setminus 
{\mathcal U}_{\delta}(x)}$, 
we can change the contour of integrals (\ref{the integral in zeta_1}) 
to the curve $\Gamma_{x, z'}$ defined by (\ref{change of variables}) 
since $\sin\theta \le \sin \theta_0=a_0$.  
For simplicity we write $\sigma_1=\rho$, $\sigma_2= \zeta_2$ and 
$\sigma=(\sigma_1, \sigma_2)$, and we set 
\begin{align}
f(\sigma) &= \sqrt{1+\sigma_1^2}\sqrt{1+\sigma_2^2},
\quad
F_k(\sigma, x, z') = Q_k(\zeta_1(\sigma_1, x, z'), \sigma_2)
\frac{1}{\sqrt{1+\sigma_1^2}}.
\label{definitions of f and F_k}
\end{align}
Then, as in the same way as section \ref{The shortest time and asymptotics} 
of \cite{transmission No1}, by (\ref{For asymptotics of the refracted part 1-0}), 
(\ref{For asymptotics of the refracted part 1-1}) and (\ref{s of I}) 
we obtain  
\begin{align}
E^{\gamma_-}_{\tau}(x, z')
&= \frac{\tau}{2(2\pi)^2\gamma_-^{3/2}}
\int_{\R^2}{e}^{-{\tilde\tau}\v{x - \tilde{z}'}f(\sigma)}
F_0(\sigma, x, z')d\sigma, 
\label{For asymptotics of the refracted part 1'}
\\
\partial_{x_k}E^{\gamma_-}_{\tau}(x, z') &= \frac{\tau^2}{2(2\pi)^2\gamma_-^{2}}
\int_{\R^2}{e}^{-{\tilde\tau}\v{x - \tilde{z}'}f(\sigma)}
F_k(\sigma, x, z')d\sigma\displaystyle\Theta_k(x, z')
\quad(k = 1, 2, 3).
\label{For asymptotics of the refracted part 2'-x'}
\end{align}
We put $\tilde{P}(\sigma, x, z') = P(\zeta_1(\sigma_1, x, z'), \sigma_2)$, 
then 
\begin{align*}
\tilde{P}(\sigma, x, z')
&=\sqrt{a_0^2-\sin^2 \theta +\sigma_1^2 \cos 2\theta + 
\frac{(1-a_0^2)\sigma_2^2}{1+\sigma_2^2}+i\sigma_1
\sqrt{1+\sigma_1^2}\sin 2\theta }.
\end{align*}
We should note that 
$F_k(\sigma, x, z')$ is continuous in $ \sigma \in \R^2$ 
and there exists a constant $C_k>0$ such that 
\begin{align}
\v{F_k(\sigma, x, z')} \leq C_k(1+\v{\sigma})^3,
\qquad(\sigma \in \R^2), 
\label{estimates of F_k including near glancing case}
\end{align}
but $F_k(\sigma, x, z')$ is not $C^{\infty}$ near $\sigma=(0,0)$ 
when $z' \in \overline{{\mathcal U}_{1}(x)\setminus 
{\mathcal U}_{\delta}(x)}$ 
because of $\tilde{P}(\sigma, x, z')$. 
For small $\vert \sigma \vert$ we will show the following continuity at $\sigma =0$:
\begin{align}
\v{F_k(\sigma, x, z') - &F_k(0, x, z')} 
\leq C(\sqrt{\v{\sigma_1}}
+\v{\sigma_2})
\nonumber
\\&\quad(\sigma, \in \R^2, \v{\sigma} \leq 2, 
z' \in \overline{{\mathcal U}_{1}(x)\setminus 
{\mathcal U}_{\delta}(x)} ).
\label{estimates of F_k including near glancing case near the origin}
\end{align}
To obtain 
(\ref{estimates of F_k including near glancing case near the origin}), 
it is enough to show 
\begin{align}
&\v{\tilde{P}(0, x, z')-\tilde{P}(\sigma, x, z')} 
\leq C(\sqrt{\v{\sigma_1}}+\v{\sigma_2}),
\label{Horder estimate of tilde{P}-1}
\\
&\Big\vert
\frac{1}{\tilde{P}(0, x, z')+a_0^2}
-\frac{1}{\tilde{P}(\sigma, x, z')+a_0^2\sqrt{1 + \sigma_1^2}}
\Big\vert
\leq C(\sqrt{\v{\sigma_1}}+\v{\sigma_2}) 
\label{Horder estimate of tilde{P}-2}
\end{align}
for $\v{\sigma} \leq 2$ and 
$z' \in \overline{{\mathcal U}_{1}(x)\setminus {\mathcal U}_{\delta}(x)}$ 
because of the definition of $Q_k$. 
Estimate 
(\ref{Horder estimate of tilde{P}-2}) follows from 
(\ref{Horder estimate of tilde{P}-1}), since 
\begin{align*}
&\hskip6mm\Big\vert
\frac{1}{\tilde{P}(0, x, z')+a_0^2}
-\frac{1}{\tilde{P}(\sigma, x, z')+a_0^2\sqrt{1 + \sigma_1^2}}
\Big\vert \\
&\leq 
\frac{ \v{ \tilde{P}(\sigma,  x, z')
-\tilde{P}(0, x, z')} 
+a_0^2\v{ 1-\sqrt{1 + \sigma_1^2}} }
{\sqrt{ \left({\rm Re}[\tilde{P}(\sigma,  x, z')]+a_0^2\sqrt{1 + \sigma_1^2} \right)^2 }
\sqrt{ \left( {\rm Re}[\tilde{P}(0, x, z')]+a_0^2 \right)^2} }
\\&
\leq \frac{\v{\tilde{P}(\sigma,  x, z')-P(0,  x, z')}
+a_0^2\v{\sigma_1}}{a_0^4}
\quad(\v{\sigma} \leq 2).
\end{align*}
Here we used the fact that ${\rm Re}[\tilde{P}(\sigma, x, z')] \geq 0$, 
which follows from the definition 
$\sqrt{X} = \v{X}^{1/2}e^{i\arg X/2}$ ($\v{\arg X} < \pi$). 
Now we shall show 
(\ref{Horder estimate of tilde{P}-1}).  
Here we consider 
\begin{align}
&\hskip6mm\tilde{P}(0, x, z')-\tilde{P}(\sigma, x, z') 
\nonumber \\
&=\tilde{P}(0, x, z')-\tilde{P}(0,\sigma_2, x, z')
+\tilde{P}(0, \sigma_2, x, z') -\tilde{P}(\sigma, x, z')
\nonumber
\\
&=\sigma_2 \int_0^1 \partial_{\sigma_2}\tilde{P}(0,t\sigma_2, x, z')\,dt
+\sigma_1 \int_0^1 \partial_{\sigma_1} 
\tilde{P}(t\sigma_1,\sigma_2, x, z')\, dt.
\label{diff of tildeP at 0} 
\end{align}
We know that 
\begin{align}
\partial_{\sigma_1}\tilde{P}(\sigma, x, z')
&=\frac{1}{2\tilde{P}(\sigma, x, z')}
\Big\{ 2\sigma_1 \cos 2\theta +i \sin 2\theta 
\Big( \frac{\sigma_1^2}{\sqrt{1+\sigma_1^2}}+\sqrt{1+\sigma_1^2} \Big)
\Big\},
\label{partial sigma1 tildeP} \\
\tilde{P}(0,\sigma_2, x, z')
&
=\sqrt{\frac{a_0^2-1}{1+\sigma_2^2}+1-\sin^2 \theta },
\nonumber
\\
\partial_{\sigma_2}\tilde{P}(0,\sigma_2, x, z')
&=\frac{(1-a_0^2)\sigma_2}{(1+\sigma_2^2)^{3/2} 
\sqrt{1-\sin^2 \theta}
\sqrt{\sigma_2^2-s(\theta)} },
\label{partial sigma2 tildeP}
\end{align}
where $s(\theta)=(\sin^2\theta-a_0^2)/(1-\sin^2\theta)$.  
To show (\ref{Horder estimate of tilde{P}-1}) by using (\ref{diff of tildeP at 0}), 
we consider 
\begin{align}
\v{\tilde{P}(\sigma,  x, z')}^4 
= \Big(a_0^2-\sin^2\theta+\sigma_1^2\cos(2\theta)
+\frac{(1-a_0^2)\sigma_2^2}{1+\sigma_2^2}
\Big)^2+\sigma_1^2(1+\sigma_1^2)\sin^2(2\theta). 
\label{absolute v. of tildeP ^4}
\end{align}
We know that 
there exists a $\epsilon>0$ such that 
$$\v{\tilde{P}(\sigma, x, z')}^4 
\ge \sigma_1^2 \sin^2(2\theta)
\ge \epsilon \sigma_1^2, 
$$
since $0 < 2\sin^{-1}(a_0\delta) \le 2\theta \le 2\theta_0 < \pi$. 
Then, it follows that there exists a constant $C$ such that 
\begin{align}
\frac{1}{\v{\tilde{P}(\sigma, x, z')}} \le \frac{C}{\v{\sigma_1}^{1/2}}
\qquad \text{ ($\v{\sigma} \le 2$, 
$z' \in \overline{{\mathcal U}_{1}(x)\setminus 
{\mathcal U}_{\delta}(x)}$).}
\label{Est of inv of tildeP for Prop2.2}
\end{align}
From 
(\ref{partial sigma1 tildeP}) and (\ref{Est of inv of tildeP for Prop2.2}) 
it follows that 
\begin{align}
\left\vert \sigma_1 \int_0^1 \partial_{\sigma_1} 
\tilde{P}(t\sigma_1,\sigma_2, x, z')dt \right\vert 
&\le C \int_0^1\frac{\v{\sigma_1}^{1/2}} 
{\sqrt{\v{t}}}\, dt 
\le 2C\v{\sigma_1}^{1/2}.
\label{est of partial sigma1 tildeP}
\end{align}
From (\ref{partial sigma2 tildeP}) and (\ref{Est of inv of tildeP for Prop2.2}) 
it follows that 
\begin{align}
\left \vert \sigma_2 \int_0^1 
\partial_{\sigma_2}\tilde{P}(0,t\sigma_2, x, z')\,dt \right \vert
&=\frac{1-a_0^2}{\sqrt{1-\sin^2 \theta}}
\int_0^1 \frac{t\sigma_2^2}{(1+t^2\sigma_2^2)^{3/2}
(t^2\sigma_2^2+\v{s(\theta)})^{1/2}}\, dt 
\nonumber\\
&\le \frac{1-a_0^2}{2\sqrt{1-\sin^2 \theta}}
\int_0^{\sigma_2^2}\frac{d\tau}{(1+\tau)^{3/2}
(\tau+\v{s(\theta)})^{1/2}}
\label{est of partial sigma2 tildeP}
\\
&
\le  \frac{1-a_0^2}{2\sqrt{1-\sin^2 \theta}}
\int_0^{\sigma_2^2} \tau^{-1/2}\, d\tau 
\nonumber\\
&\le C\v{\sigma_2}
\nonumber
\end{align}
for $\sin\theta \le a_0 <1$.  
If we apply (\ref{est of partial sigma1 tildeP}) and 
(\ref{est of partial sigma2 tildeP}) to (\ref{diff of tildeP at 0}), 
we obtain (\ref{Horder estimate of tilde{P}-1}).  
\par
Now we have prepared to show Proposition 
\ref{Asymptotics for near the total reflection angle case 1}.  
To estimate (\ref{For asymptotics of the refracted part 1'}) and  
(\ref{For asymptotics of the refracted part 2'-x'}), 
let us choose a function $\psi(\sigma)$ such that 
$\psi \in C^\infty_0(\R^2)$ with $0 \leq \psi \leq 1$, 
$\psi(\sigma) = 1$ $(\v{\sigma} \leq 1)$ and 
$\psi(\sigma) = 0$ $(\v{\sigma} \geq 3/2)$ and set 
\begin{align}
\int_{\R^2}{e}^{-{\tilde\tau}\v{x - \tilde{z}'}f(\sigma)}
F_k(\sigma, x, z')d\sigma
&= \int_{\R^2}{e}^{-{\tilde\tau}\v{x - \tilde{z}'}f(\sigma)}
F_k(\sigma, x, z')\psi(\sigma)d\sigma
\nonumber
\\&\,\,\,
+ \int_{\R^2}{e}^{-{\tilde\tau}\v{x - \tilde{z}'}f(\sigma)}
F_k(\sigma, x, z')(1-\psi(\sigma))d\sigma,  
\label{asymp est of E}
\end{align}
here $f(\sigma)$ and $F_k$ are defined by (\ref{definitions of f and F_k}).  
Since $f(\sigma) \geq 1+\v{\sigma}/4$ for $\v{\sigma} \geq 1$, it follows that
$f(\sigma) \ge 9/8+\v{\sigma}/8$ for $\v{\sigma} \ge 1$.  From this estimate and 
(\ref{estimates of F_k including near glancing case}), we have 
the estimate of the second integral of (\ref{asymp est of E}) as 
\begin{align*}
\Big\vert\int_{\R^2}e^{-{\tilde\tau}\v{x - \tilde{z}'}f(\sigma)}
F_k(\sigma, x, z')(1-\psi(\sigma))d\sigma\Big\vert
&\leq Ce^{-9{\tilde\tau\v{x - \tilde{z}'}}/8}
\int_{\R^2}(1+\v{\sigma})^3
e^{-({\tilde\tau}\v{x - \tilde{z}'}/8)\v{\sigma}}d\sigma
\nonumber\\
&\leq \frac{C_Ne^{-{\tilde\tau\v{x - \tilde{z}'}}}}
{(\tilde\tau\v{x - \tilde{z}'})^{N}}.
\end{align*}
For the first integral of (\ref{asymp est of E}) if we use Laplace method and estimate 
(\ref{estimates of F_k including near glancing case near the origin}), 
we have 
\begin{align*}
&\hskip16pt
\Big\vert\int_{\R^2}e^{-{\tilde\tau}\v{x - \tilde{z}'}f(\sigma)}
F_k(\sigma,  x, z')\psi(\sigma)d\sigma
- e^{-{\tilde\tau\v{x - \tilde{z}'}}}\left(\frac{2\pi}
{{\tilde\tau}\v{x - \tilde{z}'}}\right)
F_{k, 0}(x - \tilde{z}')
\Big\vert
\nonumber\\&
\leq C(\tilde\tau\v{x - \tilde{z}'})^{-5/4}e^{-\tilde\tau\v{x - \tilde{z}'}}. 
\end{align*}
Thus we complete the proof of 
Proposition \ref{Asymptotics for near the total reflection angle case 1}. 
\hfill
$\blacksquare$
\par
\vskip1pc
\par
\begin{figure}[h]
\begin{center}
\unitlength 0.1in
\begin{picture}( 43.6900, 22.3000)( 13.7000,-26.1000)
\put(36.3300,-3.8000){\makebox(0,0)[rt]{${\rm Im}\zeta_1$}}%
\put(59.6700,-18.3900){\makebox(0,0)[rt]{${\rm Re}\zeta_1$}}%
%
\special{pn 8}%
\special{pa 3730 2320}%
\special{pa 3730 484}%
\special{fp}%
\special{sh 1}%
\special{pa 3730 484}%
\special{pa 3710 550}%
\special{pa 3730 536}%
\special{pa 3750 550}%
\special{pa 3730 484}%
\special{fp}%
%
\special{pn 8}%
\special{pa 3690 1330}%
\special{pa 3690 1474}%
\special{fp}%
%
\special{pn 8}%
\special{pa 3770 1334}%
\special{pa 3770 1478}%
\special{fp}%
%
\special{pn 8}%
\special{ar 3730 1478 40 30  6.2831853 6.2831853}%
\special{ar 3730 1478 40 30  0.0000000 3.1415927}%
\put(36.2000,-13.5000){\makebox(0,0)[rt]{$\Gamma_{-, \varepsilon}$}}%
\put(38.9600,-15.1700){\makebox(0,0)[lb]{$\Gamma_{+, \varepsilon}$}}%
\special{pn 8}%
\special{pn 8}%
\special{pa 3760 1330}%
\special{pa 3770 1330}%
\special{pa 3780 1330}%
\special{pa 3792 1330}%
\special{pa 3802 1330}%
\special{pa 3812 1328}%
\special{pa 3824 1328}%
\special{pa 3834 1328}%
\special{pa 3844 1328}%
\special{pa 3854 1326}%
\special{pa 3866 1326}%
\special{pa 3876 1324}%
\special{pa 3886 1324}%
\special{pa 3898 1322}%
\special{pa 3908 1320}%
\special{pa 3918 1320}%
\special{pa 3928 1318}%
\special{pa 3940 1316}%
\special{pa 3950 1314}%
\special{pa 3960 1314}%
\special{pa 3972 1312}%
\special{pa 3982 1310}%
\special{pa 3992 1308}%
\special{pa 4002 1306}%
\special{pa 4014 1304}%
\special{pa 4024 1302}%
\special{pa 4034 1300}%
\special{pa 4046 1296}%
\special{pa 4056 1294}%
\special{pa 4066 1292}%
\special{pa 4078 1290}%
\special{pa 4088 1288}%
\special{pa 4098 1284}%
\special{pa 4108 1282}%
\special{pa 4120 1280}%
\special{pa 4130 1276}%
\special{pa 4140 1274}%
\special{pa 4152 1270}%
\special{pa 4162 1268}%
\special{pa 4172 1264}%
\special{pa 4182 1262}%
\special{pa 4194 1258}%
\special{pa 4204 1256}%
\special{pa 4214 1252}%
\special{pa 4226 1250}%
\special{pa 4236 1246}%
\special{pa 4246 1242}%
\special{pa 4256 1240}%
\special{pa 4268 1236}%
\special{pa 4278 1232}%
\special{pa 4288 1230}%
\special{pa 4300 1226}%
\special{pa 4310 1222}%
\special{pa 4320 1218}%
\special{pa 4330 1216}%
\special{pa 4342 1212}%
\special{pa 4352 1208}%
\special{pa 4362 1204}%
\special{pa 4374 1200}%
\special{pa 4384 1198}%
\special{pa 4394 1194}%
\special{pa 4404 1190}%
\special{pa 4416 1186}%
\special{pa 4426 1182}%
\special{pa 4436 1178}%
\special{pa 4448 1174}%
\special{pa 4458 1170}%
\special{pa 4468 1166}%
\special{pa 4480 1162}%
\special{pa 4490 1158}%
\special{pa 4500 1154}%
\special{pa 4510 1150}%
\special{pa 4522 1146}%
\special{pa 4532 1142}%
\special{pa 4542 1138}%
\special{pa 4554 1134}%
\special{pa 4564 1130}%
\special{pa 4574 1126}%
\special{pa 4584 1122}%
\special{pa 4596 1118}%
\special{pa 4606 1114}%
\special{pa 4616 1110}%
\special{pa 4628 1106}%
\special{pa 4638 1102}%
\special{pa 4648 1098}%
\special{pa 4658 1094}%
\special{pa 4670 1090}%
\special{pa 4680 1086}%
\special{pa 4690 1080}%
\special{pa 4702 1076}%
\special{pa 4712 1072}%
\special{pa 4722 1068}%
\special{pa 4732 1064}%
\special{pa 4744 1060}%
\special{pa 4754 1056}%
\special{pa 4764 1050}%
\special{pa 4776 1046}%
\special{pa 4786 1042}%
\special{pa 4796 1038}%
\special{pa 4808 1034}%
\special{pa 4818 1030}%
\special{pa 4828 1024}%
\special{pa 4838 1020}%
\special{pa 4850 1016}%
\special{pa 4860 1012}%
\special{pa 4870 1008}%
\special{pa 4882 1004}%
\special{pa 4892 998}%
\special{pa 4902 994}%
\special{pa 4912 990}%
\special{pa 4924 986}%
\special{pa 4934 980}%
\special{pa 4944 976}%
\special{pa 4956 972}%
\special{pa 4966 968}%
\special{pa 4976 964}%
\special{pa 4986 958}%
\special{pa 4998 954}%
\special{pa 5008 950}%
\special{pa 5018 946}%
\special{pa 5030 940}%
\special{pa 5040 936}%
\special{pa 5050 932}%
\special{pa 5060 928}%
\special{pa 5072 922}%
\special{pa 5082 918}%
\special{pa 5092 914}%
\special{pa 5104 910}%
\special{pa 5114 904}%
\special{pa 5124 900}%
\special{pa 5136 896}%
\special{pa 5146 892}%
\special{pa 5156 886}%
\special{pa 5166 882}%
\special{pa 5178 878}%
\special{pa 5188 872}%
\special{pa 5198 868}%
\special{pa 5210 864}%
\special{pa 5220 860}%
\special{pa 5230 854}%
\special{pa 5240 850}%
\special{pa 5252 846}%
\special{pa 5262 840}%
\special{pa 5272 836}%
\special{pa 5284 832}%
\special{pa 5294 828}%
\special{pa 5304 822}%
\special{pa 5314 818}%
\special{pa 5326 814}%
\special{pa 5336 808}%
\special{pa 5346 804}%
\special{pa 5358 800}%
\special{pa 5368 794}%
\special{pa 5378 790}%
\special{pa 5388 786}%
\special{pa 5400 780}%
\special{pa 5410 776}%
\special{pa 5420 772}%
\special{pa 5432 766}%
\special{pa 5442 762}%
\special{pa 5452 758}%
\special{pa 5462 754}%
\special{pa 5474 748}%
\special{pa 5484 744}%
\special{pa 5494 740}%
\special{pa 5506 734}%
\special{pa 5516 730}%
\special{pa 5526 726}%
\special{pa 5538 720}%
\special{pa 5548 716}%
\special{pa 5558 712}%
\special{pa 5568 706}%
\special{pa 5580 702}%
\special{pa 5590 698}%
\special{pa 5600 692}%
\special{pa 5612 688}%
\special{pa 5622 684}%
\special{pa 5632 678}%
\special{pa 5642 674}%
\special{pa 5654 668}%
\special{pa 5664 664}%
\special{pa 5674 660}%
\special{pa 5686 654}%
\special{pa 5696 650}%
\special{pa 5696 650}%
\special{pa 5696 650}%
\special{sp}%
\special{pn 8}%
\special{pa 3682 1324}%
\special{pa 3672 1324}%
\special{pa 3662 1324}%
\special{pa 3650 1324}%
\special{pa 3640 1324}%
\special{pa 3630 1322}%
\special{pa 3618 1322}%
\special{pa 3608 1322}%
\special{pa 3598 1320}%
\special{pa 3586 1320}%
\special{pa 3576 1320}%
\special{pa 3566 1318}%
\special{pa 3556 1316}%
\special{pa 3544 1316}%
\special{pa 3534 1314}%
\special{pa 3524 1314}%
\special{pa 3512 1312}%
\special{pa 3502 1310}%
\special{pa 3492 1308}%
\special{pa 3480 1306}%
\special{pa 3470 1306}%
\special{pa 3460 1304}%
\special{pa 3450 1302}%
\special{pa 3438 1300}%
\special{pa 3428 1298}%
\special{pa 3418 1296}%
\special{pa 3406 1292}%
\special{pa 3396 1290}%
\special{pa 3386 1288}%
\special{pa 3374 1286}%
\special{pa 3364 1284}%
\special{pa 3354 1280}%
\special{pa 3344 1278}%
\special{pa 3332 1276}%
\special{pa 3322 1274}%
\special{pa 3312 1270}%
\special{pa 3300 1268}%
\special{pa 3290 1264}%
\special{pa 3280 1262}%
\special{pa 3268 1258}%
\special{pa 3258 1256}%
\special{pa 3248 1252}%
\special{pa 3238 1250}%
\special{pa 3226 1246}%
\special{pa 3216 1244}%
\special{pa 3206 1240}%
\special{pa 3194 1236}%
\special{pa 3184 1234}%
\special{pa 3174 1230}%
\special{pa 3162 1226}%
\special{pa 3152 1224}%
\special{pa 3142 1220}%
\special{pa 3132 1216}%
\special{pa 3120 1212}%
\special{pa 3110 1210}%
\special{pa 3100 1206}%
\special{pa 3088 1202}%
\special{pa 3078 1198}%
\special{pa 3068 1194}%
\special{pa 3056 1192}%
\special{pa 3046 1188}%
\special{pa 3036 1184}%
\special{pa 3026 1180}%
\special{pa 3014 1176}%
\special{pa 3004 1172}%
\special{pa 2994 1168}%
\special{pa 2982 1164}%
\special{pa 2972 1160}%
\special{pa 2962 1156}%
\special{pa 2950 1152}%
\special{pa 2940 1148}%
\special{pa 2930 1146}%
\special{pa 2920 1142}%
\special{pa 2908 1138}%
\special{pa 2898 1134}%
\special{pa 2888 1128}%
\special{pa 2876 1124}%
\special{pa 2866 1120}%
\special{pa 2856 1116}%
\special{pa 2844 1112}%
\special{pa 2834 1108}%
\special{pa 2824 1104}%
\special{pa 2814 1100}%
\special{pa 2802 1096}%
\special{pa 2792 1092}%
\special{pa 2782 1088}%
\special{pa 2770 1084}%
\special{pa 2760 1080}%
\special{pa 2750 1076}%
\special{pa 2738 1072}%
\special{pa 2728 1066}%
\special{pa 2718 1062}%
\special{pa 2708 1058}%
\special{pa 2696 1054}%
\special{pa 2686 1050}%
\special{pa 2676 1046}%
\special{pa 2664 1042}%
\special{pa 2654 1036}%
\special{pa 2644 1032}%
\special{pa 2632 1028}%
\special{pa 2622 1024}%
\special{pa 2612 1020}%
\special{pa 2602 1016}%
\special{pa 2590 1010}%
\special{pa 2580 1006}%
\special{pa 2570 1002}%
\special{pa 2558 998}%
\special{pa 2548 994}%
\special{pa 2538 988}%
\special{pa 2526 984}%
\special{pa 2516 980}%
\special{pa 2506 976}%
\special{pa 2496 972}%
\special{pa 2484 966}%
\special{pa 2474 962}%
\special{pa 2464 958}%
\special{pa 2452 954}%
\special{pa 2442 948}%
\special{pa 2432 944}%
\special{pa 2420 940}%
\special{pa 2410 936}%
\special{pa 2400 932}%
\special{pa 2390 926}%
\special{pa 2378 922}%
\special{pa 2368 918}%
\special{pa 2358 914}%
\special{pa 2346 908}%
\special{pa 2336 904}%
\special{pa 2326 900}%
\special{pa 2314 894}%
\special{pa 2304 890}%
\special{pa 2294 886}%
\special{pa 2284 882}%
\special{pa 2272 876}%
\special{pa 2262 872}%
\special{pa 2252 868}%
\special{pa 2240 864}%
\special{pa 2230 858}%
\special{pa 2220 854}%
\special{pa 2208 850}%
\special{pa 2198 844}%
\special{pa 2188 840}%
\special{pa 2178 836}%
\special{pa 2166 832}%
\special{pa 2156 826}%
\special{pa 2146 822}%
\special{pa 2134 818}%
\special{pa 2124 812}%
\special{pa 2114 808}%
\special{pa 2102 804}%
\special{pa 2092 800}%
\special{pa 2082 794}%
\special{pa 2072 790}%
\special{pa 2060 786}%
\special{pa 2050 780}%
\special{pa 2040 776}%
\special{pa 2028 772}%
\special{pa 2018 766}%
\special{pa 2008 762}%
\special{pa 1996 758}%
\special{pa 1986 752}%
\special{pa 1976 748}%
\special{pa 1966 744}%
\special{pa 1954 738}%
\special{pa 1944 734}%
\special{pa 1934 730}%
\special{pa 1922 726}%
\special{pa 1912 720}%
\special{pa 1902 716}%
\special{pa 1890 712}%
\special{pa 1880 706}%
\special{pa 1870 702}%
\special{pa 1860 698}%
\special{pa 1848 692}%
\special{pa 1838 688}%
\special{pa 1828 684}%
\special{pa 1816 678}%
\special{pa 1806 674}%
\special{pa 1796 670}%
\special{pa 1784 664}%
\special{pa 1774 660}%
\special{pa 1764 656}%
\special{pa 1754 650}%
\special{pa 1754 650}%
\special{pa 1752 650}%
\special{pa 1752 650}%
\special{pa 1752 650}%
\special{pa 1752 650}%
\special{pa 1752 650}%
\special{pa 1752 650}%
\special{pa 1752 650}%
\special{pa 1752 650}%
\special{pa 1752 650}%
\special{pa 1752 650}%
\special{pa 1752 650}%
\special{pa 1752 650}%
\special{pa 1752 650}%
\special{pa 1752 650}%
\special{pa 1752 650}%
\special{pa 1752 650}%
\special{pa 1752 650}%
\special{pa 1752 650}%
\special{pa 1752 650}%
\special{pa 1750 650}%
\special{pa 1750 650}%
\special{pa 1750 650}%
\special{pa 1750 650}%
\special{pa 1750 650}%
\special{pa 1750 650}%
\special{pa 1750 650}%
\special{pa 1750 650}%
\special{pa 1750 650}%
\special{pa 1750 650}%
\special{pa 1750 648}%
\special{pa 1750 648}%
\special{pa 1750 648}%
\special{pa 1750 648}%
\special{pa 1750 648}%
\special{pa 1750 648}%
\special{pa 1750 648}%
\special{pa 1750 648}%
\special{pa 1750 648}%
\special{pa 1748 648}%
\special{pa 1748 648}%
\special{pa 1748 648}%
\special{pa 1748 648}%
\special{pa 1748 648}%
\special{pa 1748 648}%
\special{pa 1748 648}%
\special{pa 1748 648}%
\special{pa 1748 648}%
\special{pa 1748 648}%
\special{pa 1748 648}%
\special{pa 1748 648}%
\special{pa 1748 648}%
\special{pa 1748 648}%
\special{pa 1748 648}%
\special{pa 1748 648}%
\special{pa 1748 648}%
\special{pa 1748 648}%
\special{pa 1748 648}%
\special{pa 1746 648}%
\special{pa 1746 648}%
\special{pa 1746 648}%
\special{pa 1746 648}%
\special{pa 1746 648}%
\special{pa 1746 648}%
\special{pa 1746 648}%
\special{pa 1746 648}%
\special{pa 1746 648}%
\special{pa 1746 648}%
\special{pa 1746 648}%
\special{pa 1746 648}%
\special{pa 1746 648}%
\special{pa 1746 648}%
\special{pa 1746 648}%
\special{pa 1746 646}%
\special{pa 1746 646}%
\special{pa 1746 646}%
\special{pa 1746 646}%
\special{pa 1744 646}%
\special{pa 1744 646}%
\special{pa 1744 646}%
\special{pa 1744 646}%
\special{pa 1744 646}%
\special{pa 1744 646}%
\special{pa 1744 646}%
\special{pa 1744 646}%
\special{pa 1744 646}%
\special{sp}%
\put(36.7000,-16.1000){\makebox(0,0){$\Gamma_{\varepsilon}$}}%
\put(28.2500,-18.8000){\makebox(0,0){$b_0(\zeta_2)$}}%
%
\special{pn 8}%
\special{pa 4418 654}%
\special{pa 3764 1292}%
\special{dt 0.045}%
\special{sh 1}%
\special{pa 3764 1292}%
\special{pa 3826 1260}%
\special{pa 3802 1256}%
\special{pa 3798 1232}%
\special{pa 3764 1292}%
\special{fp}%
\put(44.8200,-5.9600){\makebox(0,0){$\sin\theta$}}%
%
\special{pn 8}%
\special{pa 1692 1840}%
\special{pa 1692 1840}%
\special{fp}%
%
\special{pn 8}%
\special{pa 1740 2078}%
\special{pa 5740 2078}%
\special{fp}%
\special{sh 1}%
\special{pa 5740 2078}%
\special{pa 5672 2058}%
\special{pa 5686 2078}%
\special{pa 5672 2098}%
\special{pa 5740 2078}%
\special{fp}%
%
\special{pn 8}%
\special{pa 2810 1810}%
\special{pa 3730 1460}%
\special{dt 0.045}%
\special{sh 1}%
\special{pa 3730 1460}%
\special{pa 3662 1466}%
\special{pa 3680 1480}%
\special{pa 3676 1502}%
\special{pa 3730 1460}%
\special{fp}%
\end{picture}%
\end{center}
$\quad$\vskip-15mm
\caption{Contour of the integrals}
\label{moucyounozu}
\end{figure}
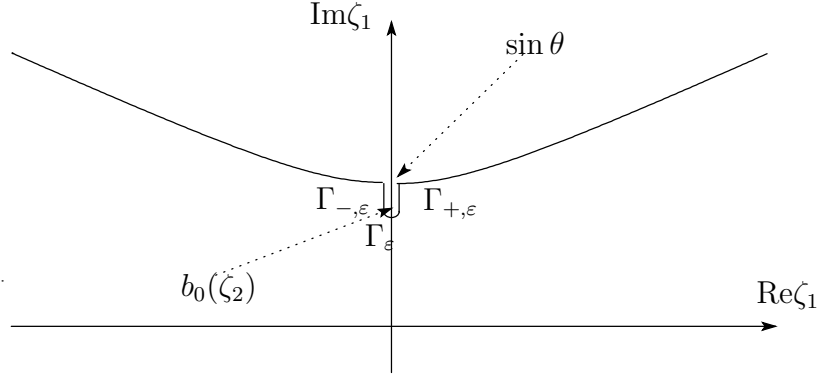

\noindent
Proof of Proposition \ref{Asymptotics for near the total reflection angle case 2}. 
Here we consider the case 
$z' \in \R^2\setminus{\mathcal U}_{1}(x)$.  
The integral $I_{\tilde\tau, k}(x-\tilde{z}', \zeta_2)$ 
in (\ref{the integral in zeta_1}) can be 
written as below: 
\begin{align}
I_{\tilde\tau, k}(x-\tilde{z}', \zeta_2)
&= \int_{\R}{e}^{-\tilde\tau r \lambda}Q_k(\zeta_1, \zeta_2)
\frac{d\zeta_1}{\sqrt{1+\zeta_1^2}},
\label{the integral in zeta_1'}
\end{align}
where $r=\vert x-\tilde{z}' \vert \sqrt{1+\zeta_2^2}$, 
$\lambda=-i(\v{x'-z'}/\vert x-\tilde{z}' \vert) \zeta_1
+(\v{x_3}/\vert x-\tilde{z}' \vert) \sqrt{1+\zeta_1^2}
=-i\sin \theta \zeta_1+\cos \theta\sqrt{1+\zeta_1^2}$ and 
$k=0, 1, 2, 3$.  
When we try to change the contour of the integrals 
$I_{\tilde{\tau}, k}(x-\tilde{z}', \zeta_2)$ ($k=0, 1, 2, 3$) 
in the same way as in the case of 
$\overline{{\mathcal U}_{1}(x)\setminus {\mathcal U}_{\delta}(x)}$, 
we need to count $a_0 < b_0(\zeta_2)$, 
and $[ib_0(\zeta_2), i\infty)$ is the branch cut of the integrands. 
Therefore, in case that $b_0(\zeta_2) < \sin \theta$, 
we consider the following contour for $\varepsilon > 0$ (see figure \ref{moucyounozu}):
\begin{align*}
&\Gamma_{\varepsilon}: \,\zeta_1=ib_0(\zeta_2)+\varepsilon e^{i\phi} 
\quad \text{($\pi \le \phi \le 2\pi$)}, 
\nonumber\\
&\Gamma_{+, \varepsilon}: \,\zeta_1=e^{i\pi/2}w+\varepsilon 
\quad \text{($b_0(\zeta_2) \le w \le \sin \theta$)},
\nonumber\\
&\Gamma_{-, \varepsilon}: \,\zeta_1=e^{i\pi/2}w-\varepsilon 
\quad \text{($\sin \theta \ge w \ge b_0(\zeta_2)$)}.
\end{align*}
%
%
%
When $\zeta_1 \in \Gamma_{\pm, \varepsilon}$, 
$\zeta_1^2=-w^2\pm 2e^{\pi{i}/2}w\varepsilon+O(\varepsilon^2)$, 
$1+\zeta_1^2=1-w^2+O(\varepsilon)$ as 
$\varepsilon \downarrow 0$.  Thus we have 
\begin{align}
P(\zeta_1, \zeta_2) &=
\sqrt{\zeta_1^2+b_0(\zeta_2)^2} \to 
\sqrt{w^2-b_0(\zeta_2)^2}e^{\pm\pi{i}/2} \quad 
\text{($\varepsilon \downarrow 0$) }
\nonumber
\end{align}
for $b_0(\zeta_2) \le w \le \sin \theta$.  
If we put 
$X_0(\zeta_2)=4\sqrt{\gamma_-}\sqrt{1+\zeta_2^2}$ and 
$X_1(\zeta_2)=4\sqrt{\gamma_-}(1+\zeta_2^2)$.  
Then we have 
\begin{align}
Q_0(\zeta_1, \zeta_2) \vert_{\zeta_1\in \Gamma_{\pm, \varepsilon} }
&=\frac{X_0(\zeta_2)\sqrt{1-w^2}\sqrt{w^2-b_0(\zeta_2)^2}e^{\pm i \pi/2}}
{\sqrt{w^2-b_0(\zeta_2)^2}e^{\pm i \pi/2}+a_0^2\sqrt{1-w^2}}
+O(\varepsilon) 
\nonumber\\
&=:Q_0^{\pm}(w,\zeta_2)+O(\varepsilon) \qquad 
\text{($\varepsilon \downarrow 0$), }
\nonumber\\
Q_k(\zeta_1, \zeta_2) \vert_{\zeta_1\in \Gamma_{\pm, \varepsilon} }
&=\frac{-X_1(\zeta_2)w
\sqrt{1-w^2}\sqrt{w^2-b_0(\zeta_2)^2}e^{\pm i \pi/2}}
{\sqrt{w^2-b_0(\zeta_2)^2}e^{\pm i \pi/2}+a_0^2\sqrt{1-w^2}}
+O(\varepsilon)
\nonumber\\
&=:Q_k^{\pm}(w,\zeta_2)+O(\varepsilon) \qquad 
\text{($\varepsilon \downarrow 0$, $k = 1, 2$)}
\nonumber\\
Q_3(\zeta_1, \zeta_2) \vert_{\zeta_1\in \Gamma_{\pm, \varepsilon} }
&=\frac{-X_1(\zeta_2)(1-w^2)\sqrt{w^2-b_0(\zeta_2)^2}e^{\pm i \pi/2}}
{\sqrt{w^2-b_0(\zeta_2)^2}e^{\pm i \pi/2}+a_0^2\sqrt{1-w^2}}
+O(\varepsilon) 
\nonumber\\
&=:Q_3^{\pm}(w,\zeta_2)+O(\varepsilon) \qquad 
\text{($\varepsilon \downarrow 0$),}
\nonumber
\\
\lambda \vert_{\zeta_1\in \Gamma_{\pm, \varepsilon} }
&=-i(\sin\theta)\zeta_1+(\cos \theta)\sqrt{1+\zeta_1^2}  
\vert_{\zeta_1\in \Gamma_{\pm, \varepsilon} }
\nonumber\\
&= (\sin \theta)w + (\cos \theta)\sqrt{1-w^2} +O(\varepsilon) 
\nonumber\\
&=: \lambda_0(w)+O(\varepsilon) 
\qquad 
\text{($\varepsilon \downarrow 0$)}.
\nonumber
\end{align}
We define the following integrals 
$I_{\tilde\tau, k}^s(x-\tilde{z}', \zeta_2)$ and 
$I_{\tilde\tau, k}^m(x-\tilde{z}', \zeta_2)$ for $k=0, 1, 2, 3$ 
as below:
\begin{align}
I_{\tilde\tau, k}^s(x-\tilde{z}', \zeta_2)
&=\int_{\Gamma_{x, z'}}e^{-\tilde{\tau}r\lambda}
Q_k(\zeta_1, \zeta_2)\frac{d\zeta_1}{\sqrt{1+\zeta_1^2}},
\label{sug}\\
I^m_{\tilde{\tau},k}(x-\tilde{z}',\zeta_2)
&=\lim_{\varepsilon \downarrow 0} 
\int_{\Gamma_{+, \varepsilon}\cup\Gamma_{-, \varepsilon}
\cup \Gamma_{\varepsilon}}
\frac{Q_k(\zeta_1, \zeta_2)}{\sqrt{1+\zeta_1^2}} e^{-\tilde{\tau}r\lambda}\,
d\zeta_1 \nonumber\\
&=i\int_{b_0(\zeta_1)}^{\sin\theta} 
\frac{Q_k^+(w,\zeta_2)-Q_k^-(w,\zeta_2)}{\sqrt{1-w^2}}e^{-\tilde{\tau}r\lambda_0(w)}\,dw.
\label{appe}
\end{align}
Then we can change the contour of (\ref{the integral in zeta_1'}) as 
\begin{align*}
I_{\tilde\tau, k}(x-\tilde{z}', \zeta_2)
=I^s_{\tilde{\tau},k}(x-\tilde{z}',\zeta_2)
+I^m_{\tilde{\tau},k}(x-\tilde{z}',\zeta_2).
\end{align*}
Thus (\ref{For asymptotics of the refracted part 1-0}) and  
(\ref{For asymptotics of the refracted part 1-1}) are reduced to 
\begin{align*}
{E}_{\tau}^{\gamma_-}(x, z') 
&= {E}_{\tau, 0}^{s,\gamma_-}(x, z') +{E}_{\tau, 0}^{m,\gamma_-}(x, z'),
\\
\partial_{x_k}{E}_{\tau}^{\gamma_-}(x, z') 
&= {E}_{\tau, k}^{s,\gamma_-}(x, z') +{E}_{\tau, k}^{m,\gamma_-}(x, z')
\qquad(k = 1, 2, 3), 
\end{align*}
where
\begin{align}
{E}_{\tau, 0}^{\alpha,\gamma_-}(x, z') &= 
\frac{\tau}{2(2\pi)^2\gamma_-^{3/2}}\int_{\R}
I^\alpha_{\tilde\tau, 0}(x-\tilde{z}', \zeta_2)d\zeta_2 \quad(\alpha = s, m),
\label{E^alpha_tau-0}
\\
{E}_{\tau, k}^{\alpha,\gamma_-}(x, z') &= 
\frac{\tau^2}{2(2\pi)^2\gamma_-^{2}}\int_{\R}
I^\alpha_{\tilde\tau, k}(x-\tilde{z}', \zeta_2)d\zeta_2\Theta_k(x, z') 
\quad(\alpha = s, m, k = 1, 2, 3).
\label{E^alpha_tau-1}
\end{align}
\par

At first, we shall consider ${E}_{\tau, k}^{s,\gamma_-}(x, z')$. From (\ref{sug}), 
they are reduced to similar 
forms to (\ref{For asymptotics of the refracted part 1'}) and 
(\ref{For asymptotics of the refracted part 2'-x'}).
Hence, if we prove estimates corresponding to (\ref{Horder estimate of tilde{P}-1}) for 
$z' \in \R^2\setminus{\mathcal U}_{1}(x)$, the same argument as for 
(\ref{For asymptotics of the refracted part 1'}) and 
(\ref{For asymptotics of the refracted part 2'-x'}) works. 
Thus we should show
\begin{align}
\frac{1}{\v{\tilde{P}(\sigma, x, z')}} \le \frac{C}{\v{\sigma_1}^{1/2}}
\qquad \text{ ($\v{\sigma} \le 2$, 
$z' \in \R^2\setminus{\mathcal U}_{1}(x)$).}
\label{Est of inv of tildeP for Prop2.3}
\end{align}

\par

If we set 
$Y=a_0^2-\sin^2\theta+\sigma_1^2\cos(2\theta)+\frac{(1-a_0^2)\sigma_2^2}{1+\sigma_2^2}$
and $Y_0=\sigma_1^2+\frac{1-a_0^2}{1+\sigma_2^2}$, then 
$Y=-Y_0 + \cos^2\theta(1+2\sigma_1^2)$ and 
$0 < (1-a_0^2)/5 \le Y_0 \le 5$ as $\v{\sigma} \le 2$.  
From 
$$
\v{\tilde{P}(\sigma, x, z')}^4 \ge 
Y^2 \ge Y_0^2- 2 Y_0\cos^2\theta(1+2\sigma_1^2) \ge Y_0(Y_0-18\cos^2\theta), 
$$
it follows that 
$\v{\tilde{P}(\sigma, x, z')}^4 \ge (1-a_0^2)^2/50$ for $\cos^2 \theta \le (1-a_0^2)/180$. 
When $\cos^2 \theta > (1-a_0^2)/180$, 
there exists a $\epsilon > 0$ such that 
$\v{\tilde{P}(\sigma, x, z')}^4 \ge \epsilon \sigma_1^2$, 
since we know that
$\v{\tilde{P}(\sigma, x, z')}^4 \ge \sigma_1^2\sin^2(2\theta)$ 
from (\ref{absolute v. of tildeP ^4}).  
Thus, we obtain (\ref{Est of inv of tildeP for Prop2.3}).

\par

By using (\ref{diff of tildeP at 0}), we shall estimate 
$\tilde{P}(0, x, z')-\tilde{P}(\sigma, x, z')$.  
Then, the first term is 
\begin{align}
\sigma_2 \int_0^1 \partial_{\sigma_2}\tilde{P}(0,t\sigma_2, x, z')\,dt
&=\int_0^1 
\frac{(1-a_0^2)t\sigma_2^2}{(1+t^2\sigma_2^2)^{3/2} 
\sqrt{1-\sin^2 \theta}
\sqrt{t^2\sigma_2^2-s(\theta)} }\,dt
\nonumber
\\
&=\frac{1-a_0^2}{2\sqrt{1-\sin^2 \theta}}
\int_0^{\sigma_2^2} 
\frac{d\tau}{(1+\tau)^{3/2}(\tau-s(\theta))^{1/2}}.
\nonumber
\end{align}
Thus we have 
\begin{align}
\left\vert \sigma_2 \int_0^1 
\partial_{\sigma_2}\tilde{P}(0,t\sigma_2, x, z')\,dt \right\vert
\le \frac{1-a_0^2}{2\sqrt{1-\sin^2 \theta}}
\int_0^{\sigma_2^2} \v{ \tau-s(\theta)}^{-1/2}\,d\tau. 
\nonumber
\end{align}
When $s(\theta) < \sigma_2^2$, we have
\begin{align}
&\hskip6mm \left\vert \sigma_2 \int_0^1 
\partial_{\sigma_2}\tilde{P}(0,t\sigma_2, x, z')\,dt \right\vert
\nonumber\\
&\le \frac{1-a_0^2}{2\sqrt{1-\sin^2 \theta}}
\Big( 
\int_0^{s(\theta)} (s(\theta)-\tau)^{-1/2}\,d\tau 
+\int_{s(\theta)}^{\sigma_2^2} 
(\tau-s(\theta))^{-1/2}\,d\tau \Big)
\nonumber\\
&=\frac{1-a_0^2}{2\sqrt{1-\sin^2 \theta}}
\Big(
2s(\theta)^{1/2}+2(\sigma_2^2-s(\theta))^{1/2}\Big)
\nonumber\\
&\le C(s(\theta)^{1/2}+\v{\sigma_2})
\nonumber\\
&\le C \v{\sigma_2}. 
\nonumber
\end{align}
When $s(\theta) \ge \sigma_2^2$, we have
\begin{align}
\left\vert \sigma_2 \int_0^1 
\partial_{\sigma_2}\tilde{P}(0,t\sigma_2, x, z')\,dt \right\vert
&\le \frac{1-a_0^2}{2\sqrt{1-\sin^2 \theta}}
\Big( 
\int_0^{\sigma_2^2} (s(\theta)-\tau)^{-1/2}\,d\tau 
\Big) 
\nonumber\\
&\le \frac{1-a_0^2}{\sqrt{1-\sin^2 \theta}}
\Big\{ \sqrt{s(\theta)}-\sqrt{s(\theta)-\sigma_2^2} \Big\}
\nonumber\\
&\le C \frac{(\sigma_2/\sqrt{s(\theta)})\sigma_2}
{1+\sqrt{1-\sigma_2^2/s(\theta)}}
\nonumber\\
&\le C \v{\sigma_2}. 
\nonumber
\end{align}
Next, the second term can be estimated as follows 
\begin{align}
&\hskip6mm\v{\sigma_1 \partial_{\sigma_1}\tilde{P}(t\sigma_1, \sigma_2, x, z')}^2
\nonumber\\
&=\frac{\sigma_1^2}{4\v{\tilde{P}(t\sigma_1,\sigma_2, x, z')}^2}
\Big\{ 4(t\sigma_1)^2 \cos^2(2\theta)+\sin^2(2\theta)
\Big(\frac{(t\sigma_1)^2}{\sqrt{1+(t\sigma_1)^2}}
+ \sqrt{1+(t\sigma_1)^2}\Big)^2 \Big\}
\nonumber\\
&\le \frac{C\sigma_1^2}{4\v{t\sigma_1}}
\nonumber\\&
\le C\frac{\v{\sigma_1}}{\v{t}},
\nonumber
\end{align}
by using (\ref{Est of inv of tildeP for Prop2.3}).  
Then, 
\begin{align}
\left\vert \sigma_1 \int_0^1 
\partial_{\sigma_1}\tilde{P}(t\sigma_1, \sigma_2, x, z')\,dt \right\vert 
\le \tilde{C}\v{\sigma_1}^{1/2} \int_0^1 \frac{1}{\sqrt{t}}\,dt 
=2\tilde{C}\v{\sigma_1}^{1/2}.
\nonumber
\end{align}
Thus we have (\ref{Horder estimate of tilde{P}-1}) for $\v{\sigma} \leq 2$ and 
$z' \in \R^2\setminus{\mathcal U}_{1}(x)$.  
Using this, we can follow the argument getting the estimates for 
(\ref{For asymptotics of the refracted part 1'}) and 
(\ref{For asymptotics of the refracted part 2'-x'}).
Thus, we obtain
\begin{align}
\v{{E}_{\tau, 0}^{s,\gamma_-}(x, z')}
&\le Ce^{-\tau \frac{\v{x-\tilde{z}'}}{\sqrt{\gamma_-}}} 
\qquad \text{($x \in \overline{D}$, 
$z' \in \R^2\setminus{\mathcal U}_{1}(x)$),}
\label{Est of E sug}
\\
\v{{E}_{\tau, k}^{s,\gamma_-}(x, z')}
&\le C{\tau}e^{-\tau \frac{\v{x-\tilde{z}'}}{\sqrt{\gamma_-}}} 
\qquad \text{($x \in \overline{D}$, 
$z' \in \R^2\setminus{\mathcal U}_{1}(x)$, $k = 1, 2, 3$).}
\label{Est of E sug'}
\end{align}

\par

From now, we shall estimate ${E}_{\tau, k}^{m,\gamma_-}(x, z')$.
From (\ref{appe}), the integral 
$I^m_{\tilde{\tau}, k}(x-\tilde{z}',\zeta_2)$ ($k = 0, \ldots, 3$) 
can be expressed by 
\begin{align}
I^m_{\tilde{\tau},0}(x-\tilde{z}',\zeta_2)&=
-8\sqrt{\gamma_-}\sqrt{1+\zeta_2^2}
\int_{b_0(\zeta_2)}^{\sin\theta}
G(w, \zeta_2)
e^{-\tilde{\tau} r \lambda_0(w)}\, dw,
\nonumber
\\
I^m_{\tilde{\tau},k}(x-\tilde{z}',\zeta_2)&=
8\sqrt{\gamma_-}(1+\zeta_2^2)
\int_{b_0(\zeta_2)}^{\sin\theta}
wG(w, \zeta_2)
e^{-\tilde{\tau} r \lambda_0(w)}\, dw, \quad(k = 1, 2), 
\nonumber\\
I^m_{\tilde{\tau},3}(x-\tilde{z}',\zeta_2)&=
8\sqrt{\gamma_-}(1+\zeta_2^2)
\int_{b_0(\zeta_2)}^{\sin\theta}
\sqrt{1 -w^2}G(w, \zeta_2)
e^{-\tilde{\tau} r \lambda_0(w)}\, dw,
\nonumber
\end{align}
where 
$$
G(w, \zeta_2) = \frac{a_0^2\sqrt{1 -w^2}\sqrt{w^2-b_0(\zeta_2)^2}}
{a_0^4(1 -w^2)+ \vert w^2 -b_0(\zeta_2)^2 \vert}. 
$$
Since $0 \leq G(w, \zeta_2) \leq 1/2$, 
we have 
\begin{align}
\v{ I^m_{\tilde{\tau},0}(x-\tilde{z}',\zeta_2)}
&\le 4\sqrt{\gamma_-}\sqrt{1+\zeta_2^2}
\int_{\sin \theta_0}^{\sin\theta}
e^{-\tilde{\tau} r (w\sin\theta+\sqrt{1-w^2}\cos \theta)}\,dw, 
\nonumber 
\end{align}
where $\sin \theta_0=a_0 < b_0(\zeta_2)$ is used. 
Moreover, from the change of variable $w=\sin \alpha$ and 
the relation $w\sin\theta+\sqrt{1-w^2}\cos \theta
=\sin \alpha \sin \theta +\cos \alpha \cos \theta
=\cos(\theta- \alpha)$ it follows that 
\begin{align}
\v{I^m_{\tilde{\tau},0}(x-\tilde{z}',\zeta_2)}
\le 4\sqrt{\gamma_-}\sqrt{1+\zeta_2^2}
\int_{\theta_0}^{\theta} \cos \alpha \,
e^{-{\tau} T(\alpha)
\sqrt{1+\zeta_2^2}}\, d\alpha,
\nonumber
\end{align}
where we used expression (\ref{expression of T(alpha)}) of
$T(\alpha)=T_{x, z'}(\alpha)$ defined by (\ref{T(alpha)}).  
In the same way as above, we have 
\begin{align}
\v{I^m_{\tilde{\tau}, k}(x-\tilde{z}',\zeta_2)}
&\le 4\sqrt{\gamma_-}(1+\zeta_2^2)
\int_{\theta_0}^{\theta}\sin \alpha \cos \alpha \,
e^{-{\tau}T(\alpha)\sqrt{1+\zeta_2^2}}
\,d\alpha \quad(k = 1, 2),
\nonumber\\
\v{I^m_{\tilde{\tau},3}(x-\tilde{z}',\zeta_2)}
&\le 
4\sqrt{\gamma_-}(1+\zeta_2^2)
\int_{\theta_0}^{\theta}\cos^2\alpha \,
e^{-{\tau}T(\alpha)\sqrt{1+\zeta_2^2}}
\,d\alpha.
\nonumber
\end{align}
Since $2\theta - \theta_0 - \alpha \geq \alpha - \theta_0 \ge 0$ for 
$\theta_0 \leq \alpha \leq \theta$, noting $\sin t \geq 2t/\pi$ for $0 \leq t \leq \pi/2$, 
we obtain 
\begin{align*}
T(\alpha) - T(\theta_0) &= \frac{\v{\tilde{z}'-x}}{\sqrt{\gamma_-}}
(\cos(\theta-\alpha)-\cos(\theta-\theta_0))
\\
&= \frac{2\v{\tilde{z}'-x}}{\sqrt{\gamma_-}}\sin\Big(\frac{2\theta-\theta_0-\alpha}{2}\Big)
\sin\Big(\frac{\alpha-\theta_0}{2}\Big)
\\&
\geq \frac{2\v{\tilde{z}'-x}}{\pi^2\sqrt{\gamma_-}}(2\theta-\theta_0-\alpha)(\alpha-\theta_0)
\\
&\geq \frac{2\v{\tilde{z}'-x}}{\pi^2\sqrt{\gamma_-}}(\alpha-\theta_0)^2.
\end{align*}
Hence, we have 
\begin{align*}
\v{I^m_{\tilde{\tau},k}(x-\tilde{z}',\zeta_2)}
&\le C(1+\zeta_2^2)^{\frac{k_0+1}{2}}e^{-\tau T(\theta_0)\sqrt{1+\zeta_2^2}}
\int_{\theta_0}^{\theta}e^{-\frac{2\tau\v{\tilde{z}'-x}\sqrt{1+\zeta_2^2}}{\pi^2\sqrt{\gamma_-}}
(\alpha-\theta_0)^2}d\alpha
\\&
\leq C\tau^{-1/2}(1+\zeta_2^2)^{\frac{2k_0+1}{4}}e^{-\tau T(\theta_0)\sqrt{1+\zeta_2^2}}, 
\nonumber
\end{align*}
where  $k_0=0$ for $k=0$ and 
$k_0=1$ for $k=1, 2, 3$. 
Thus we have 
\begin{align}
\left\vert
\int_{\R} I^m_{\tilde{\tau},k}(x-\tilde{z}',\zeta_2)d\zeta_2
\right\vert
\leq \tilde{C}\tau^{-1/2}
\int_0^{\infty}e^{-\tau T(\theta_0)\sqrt{1+s^2}} 
(1+s^2)^{(2{k_0}+1)/4}\, ds. 
\nonumber
\end{align}
Since $\sqrt{1+s^2} \geq 1+s^2/3$ for $0 \leq s \leq 1$, it follows that 
\begin{align}
&\hskip6mm\int_{0}^{\infty}e^{-\tau T(\theta_0)\sqrt{1+s^2}} 
(1+s^2)^{(2{k_0}+1)/4}\, ds 
\nonumber\\
&\leq \int_0^{1}2^{(2{k_0}+1)/4} 
e^{-\tau T(\theta_0)}e^{-\tau T(\theta_0)s^2/3}ds
+\int_{1}^{\infty}(1+s^2)^{(2{k_0}-1)/4}\sqrt{2s^2} 
e^{-\tau T(\theta_0)\sqrt{1+s^2}} \,ds 
\nonumber \\
&\le \frac{2^{(2{k_0}-3)/4}\sqrt{3\pi}e^{-\tau T(\theta_0)}}{\sqrt{\tau T(\theta_0)}}
+\sqrt{2}\int_{\sqrt{2}}^{\infty} 
e^{-\tau T(\theta_0)\tilde{s}}
\tilde{s}^{{k_0}+1/2}
\,d\tilde{s} \qquad (\tilde{s}=\sqrt{1+s^2})
\nonumber \\
&\le C{\tau}^{-1/2}e^{-\tau T(\theta_0)}
+\sqrt{2}\int_{0}^{\infty} 
e^{-\tau T(\theta_0)(\sqrt{2}+s)}
(\sqrt{2}+s)^{{k_0}+1/2}
\,ds
\nonumber 
\\
&\le \tilde{C}{\tau}^{-1/2}e^{-\tau T(\theta_0)}.
\nonumber 
\end{align}
Thus we have
\begin{align}
\left\vert
\int_{\R} I^m_{\tilde{\tau},k}(x-\tilde{z}',\zeta_2) d\zeta_2
\right\vert
\le {C}{\tau}^{-1}e^{-\tau T(\theta_0)},
\nonumber
\end{align}
which means that 
\begin{align}
\v{{E}_{\tau, k}^{m,\gamma_-}(x, z') } 
\le C_k \tau^{k_0} e^{-\tau T(\theta_0)} 
\quad(k = 0, 1, 2, 3)
\label{Est of E appe}
\end{align} 
from (\ref{E^alpha_tau-0}) and (\ref{E^alpha_tau-1}). 
Since $\frac{\v{x-\tilde{z}'}}{\sqrt{\gamma_-}} \ge 
\frac{\v{x-\tilde{z}'}}{\sqrt{\gamma_-}}\cos(\theta-\theta_0)=T(\theta_0)$, 
Proposition 2.3 is proved by (\ref{Est of E sug}) - (\ref{Est of E appe}).  
\hfill$\blacksquare$
\setcounter{equation}{0}
\section{The optical distance and asymptotics of $\Phi_{\tau}(x, y)$}
\label{The shortest time and asymptotics} 
\vskip0pt\noindent

For $(x, y) \in \R^3_-\times\R^3_+$, we define $\tilde{l}_{x, y}(z')$ by 
\begin{align}
\tilde{l}_{x, y}(z') = 
\begin{cases} 
l_{x, y}(z') & (z' \in {\mathcal U}_{1}(x)), \\[2mm]
\displaystyle
\frac{\v{x_3}\cos\theta_0}{\sqrt{\gamma_-}}
+\frac{\v{x' - z'}+\v{\tilde{z}' - y}}{\sqrt{\gamma_+}} & 
(z' \in \R^2{\setminus}{\mathcal U}_{1}(x)),
\end{cases}
\label{def of tilde{l}_{x, y}(z')}
\end{align}
where $0 < \theta_0 < \pi/2$ is given by $\sin\theta_0 = a_0 < 1$.
Note that 
\begin{align}
T(\theta_0)+\frac{\v{\tilde{z}'- y}}{\sqrt{\gamma_+}} 
= \frac{1}{\sqrt{\gamma_-}}\Big(\v{x_3}\cos\theta_0 + \v{z'-x'}\sin\theta_0\Big)
+\frac{\v{\tilde{z}'- y}}{\sqrt{\gamma_+}}
= \tilde{l}_{x, y}(z')
\label{why tilde{l}_{x, y}(z') is defined as it is}
\end{align}
for $z' \in \R^2{\setminus}{\mathcal U}_{1}(x)$. 

\par

Proposition \ref{Asymptotics for near the total reflection angle case 2} shows that 
$\tilde{l}_{x, y}(z')$ gives the time in which the waves travel from 
$x$ to $y$ via $\tilde{z} \in \partial\R^3_+$ if $ z' \in \R^2\setminus{\mathcal U}_{1}(x)$.
This arrival time $\tilde{l}_{x, y}(z')$ is different from ${l}_{x, y}(z')$, which is caused
by the total reflection phenomena.

\par

Let us explain the meaning of $\tilde{l}_{x, y}(z')$ for $z' \in \R^2\setminus{\mathcal U}_{1}(x)$.
Since 
$\v{x'-z'}/\v{x-\tilde{z}'} > \sin \theta_0$, there exists a point 
$z'_0 = z'_0(x, z')\in \R^2$ on the line segment $x'z'$ such that 
$\v{x'-z_0'}/\v{x-\tilde{z_0}'} = \sin \theta_0$ and  
$\v{x'-z'}=\v{x'-z_0'}+\v{z_0'-z'}$. 
Note that $\tilde{l}_{x, y}(z')$ is written by
\begin{align}
\tilde{l}_{x, y}(z')
&=\frac{\cos\theta_0 }{\sqrt{\gamma_-}}
\frac{\v{x_3}}{\v{x-\tilde{z}'_0}}
\v{x-\tilde{z}'_0}
+\frac{\v{x'-z_0'}+\v{z_0'-z'}}{\sqrt{\gamma_+}} 
+\frac{\v{\tilde{z}'-y}}{\sqrt{\gamma_+}}
\nonumber \\
&=\frac{\cos^2\theta_0 }{\sqrt{\gamma_-}}
\v{x-\tilde{z}'_0}
+\frac{\sqrt{\gamma_-}}{\sqrt{\gamma_+}}
\frac{\v{x'-z_0'}}{\v{x-\tilde{z}'_0}}
\frac{\v{x-\tilde{z}'_0}}{\sqrt{\gamma_-}}
+\frac{\v{z_0'-z'}}{\sqrt{\gamma_+}}
+\frac{\v{\tilde{z}'-y}}{\sqrt{\gamma_+}}
\nonumber\\
&=\frac{\v{x-\tilde{z}'_0}}{\sqrt{\gamma_-}}
+\frac{\v{z_0'-z'}+\v{\tilde{z}'-y}}{\sqrt{\gamma_+}}. 
\label{another form of tilde{l}_{x, y}}
\end{align}
This means that if the total reflection is caused, i.e. 
$z' \in \R^2\setminus{\mathcal U}_{1}(x)$, the waves emanating from $x$ and 
arriving at $y$ via $\tilde{z}'$ go to $\tilde{z}_0'\in \partial\R^3_+$ first, 
move to $\tilde{z}'$ along the transmission boundary $\partial\R^3_+$, 
and travel to $y$ in $\R^3_+$.

\par

To obtain Proposition \ref{Asymptotics of the refracted part of the gradient of the FS}, 
we need to find $\inf_{z' \in \R^2}\tilde{l}_{x, y}(z')$. 
From Lemma 4.1 in \cite{transmission No1}, $l(x, y)$ in (\ref{definition of l(x, y)}) is attained
by only one point $ z'(x, y) $ which is $C^\infty$ for $(x, y) \in \R^3_-\times\R^3_+$. 
Note that this point $ z'(x, y) $ is determined by Snell's law
\begin{align}
\frac{\sin\theta_-}{\sqrt{\gamma_-}} = \frac{\sin\theta_+}{\sqrt{\gamma_+}}, 
\label{Snell's law-1'}
\end{align}
where $0 \leq \theta_\pm < \pi/2$ is taken by
\begin{align}
\sin\theta_- = \frac{\v{z'(x, y) - x'}}{\v{\tilde{z}'(x, y) - x}},
\qquad
\sin\theta_+ = \frac{\v{z'(x, y) - y'}}{\v{\tilde{z}'(x, y) - y}}. 
\label{definition of angle for refracted wave}
\end{align}
As in the proof of Lemma 4.1 in \cite{transmission No1}, 
\begin{equation}
\v{z(x, y) - x'} \leq \v{x' - y'}, \qquad \v{z(x, y) - y'} \leq \v{x' - y'}, 
\label{z(x, y) is in x'y'}
\end{equation}
since
$
l(x, y) = \inf\{ l_{x, y}(z') \mid z' \in \R^2, \v{z' - x'} \leq \v{x' - y'}, 
\v{z' - y'} \leq \v{x' - y'}\}.
$

\par

Here, we show the following properties of $z'(x, y)$ and the function $\tilde{l}_{x, y}(z')$.
\begin{Lemma}\label{the shortest length in the case}
(1) For any $x, y \in \R^3$ with $x_3 < 0$ and $y_3 > 0$,
$\inf_{z' \in \R^2}\tilde{l}_{x, y}(z') = l(x, y)$, and this infimum is 
attained at only $z' = z'(x, y)$. 
\par\noindent
(2) There exists a constant $0 < \delta_0 < 1$ such that
for any $(x, y) \in \overline{D}\times\overline{B}$, 
$z'(x, y) \in \overline{{\mathcal U}_{\delta_0}(x)}$, that is 
\begin{align}
\v{x' - z'(x, y)} \leq a_0\delta_0\v{x - \tilde{z}'(x, y)}
\qquad((x, y) \in \overline{D}\times\overline{B}), 
\label{z'(x, y) is apart from the sets giving critical angle}
\end{align}
and for any $\delta_1 > 0$ with $\delta_0 < \delta_1 $, there exists a constant $c_0 > 0$ such that
\begin{align}
\tilde{l}_{x, y}(z') \geq l(x, y)+c_0\v{z' - z'(x, y)} \qquad 
((x, y) \in \overline{D}\times\overline{B}, 
z' \in \R^2\setminus{\mathcal U}_{\delta_1}(x)).
\label{estimate of tilde{l}_{x, y}}
\end{align}
\end{Lemma}
Proof. From Lemma 4.1 in \cite{transmission No1}, 
(1) of Lemma \ref{the shortest length in the case} is obvious 
if (2) of Lemma \ref{the shortest length in the case} is obtained. 
To show (2) of Lemma \ref{the shortest length in the case}, 
we take constants $A > 0$ and $L > 0$ satisfying
\begin{align}
\v{x' - y'} \leq L, A \leq \v{x_3} \leq A^{-1}, 
A \leq \v{y_3} \leq A^{-1}
\qquad
&(x \in \overline{D}, y \in \overline{B}).
\qquad 
\label{bound of x,y}
\end{align}
From Snell's law (\ref{Snell's law-1'}), for any $x \in \overline{D}$ and $y \in \overline{B}$, 
$\theta_\pm$ defined by (\ref{definition of angle for refracted wave}) satisfies
$\sin \theta_-=\sqrt{\frac{\gamma_-}{\gamma_+}}\sin \theta_+
=\sin \theta_0 \sin \theta_+=a_0\sin \theta_+$.  
From (\ref{definition of angle for refracted wave}), (\ref{z(x, y) is in x'y'}) 
and (\ref{bound of x,y}), it follows that 
\begin{align*}
0 \le \sin \theta_+ =
\frac{\v{y'-z'}}{\sqrt{\v{y'-z'}^2+y_3^2}}
\le \frac{\v{y'-z'}}{\sqrt{\v{y'-z'}^2+A^2}}
\le \frac{L}{\sqrt{L^2+A^2}} < 1
\end{align*}
since $t \mapsto \frac{t}{\sqrt{t^2+A^2}}$ is monotone increasing for $t \geq 0$. 
Choose $0 < \theta_{max} <\pi/2$ satisfying $\sin \theta_{max}=\frac{L}{\sqrt{L^2+A^2}}$. 
Then for any $x \in \overline{D}$ and $y \in \overline{B}$, 
we have $\sin \theta_-=a_0 \sin \theta_+ 
\leq a_0 \sin \theta_{max} < a_0$. Hence, putting $\delta_0 = \sin \theta_{max}$, 
we obtain
\begin{equation}
z'(x, y) \in \overline{{\mathcal U}_{\sin\theta_{max}}(x)} 
= \overline{{\mathcal U}_{\delta_0}(x)}
\qquad
(x \in \overline{D}, y \in \overline{B}),  
\label{range of z'(x, y)}
\end{equation}
which gives (\ref{z'(x, y) is apart from the sets giving critical angle}).

\par
It suffices to show (\ref{estimate of tilde{l}_{x, y}}) for $\delta_1$ with
$\delta_0 < \delta_1 < 1$, since 
${\mathcal U}_{\delta_1}(x) \subset {\mathcal U}_{\delta_2}(x)$ for $\delta_1 < \delta_2$.
Take any $\delta_1$ with $1 > \delta_1 > \delta_0$ and put 
${\mathcal K} = \{(x, y, z') \in \overline{D}\times\overline{B}\times\R^2 \,\vert\, z' \in 
\overline{{\mathcal U}_{1}(x)\setminus{\mathcal U}_{\delta_1}(x)}
\,\}$. 
Noting (\ref{estimate for x and z' in overline{{mathcal U}_{delta}(x)}}), we obtain
\begin{align}
\frac{a_0\delta_1}{\sqrt{1-a_0^2\delta_1^2}}\v{x_3} \leq \v{x' - z'} 
\leq \frac{a_0}{\sqrt{1-a_0^2}}\v{x_3}
\qquad (z' \in \overline{{\mathcal U}_{1}(x)\setminus{\mathcal U}_{\delta_1}(x)}).
\label{estimate for x and z' in overline{{mathcal U}_{1}(x)setminus{mathcal U}_{delta_1}(x)}}
\end{align}
Compactness of $\overline{D}\times\overline{B}$, (\ref{bound of x,y}) and 
(\ref{estimate for x and z' in overline{{mathcal U}_{1}(x)setminus{mathcal U}_{delta_1}(x)}}) 
imply that ${\mathcal K}$ is compact. 
From (\ref{z'(x, y) is apart from the sets giving critical angle}), (\ref{bound of x,y}), 
(\ref{range of z'(x, y)}) and
(\ref{estimate for x and z' in overline{{mathcal U}_{1}(x)setminus{mathcal U}_{delta_1}(x)}}), 
it follows that there exists a constant $c_1 > 0$ such that
\begin{align}
\frac{l_{x, y}(z') - l(x, y)}{\v{z' - z'(x, y)}} \geq c_1
\quad ((x, y, z') \in {\mathcal K}),
\label{estimate for the l_{x, y}(z') part}
\end{align}
since the function in the above is positive and continuous on ${\mathcal K}$.

\par

Next, take an arbitrary $ (x, y, z') \in \overline{D}\times\overline{B}\times\R^2$ with
$z' \in \R^2\setminus{\mathcal U}_{1}(x)$. 
We take $z_0' = z_0'(x, z')$ as in (\ref{another form of tilde{l}_{x, y}}).
Since $y \neq \tilde{z}'$ and $z' \neq z_0'$, from
\begin{align*}
\v{y - \tilde{z}'}^2 = \Big\vert 
\v{y - \tilde{z}_0'} - (\tilde{z}'- \tilde{z}_0')
\cdot\frac{y - \tilde{z}_0'}{\v{y - \tilde{z}_0'}} \Big\vert^2 
+ \v{\tilde{z}' - \tilde{z}_0'}^2\Big(1 
- \Big(\frac{\tilde{z}' - \tilde{z}_0'}{\v{\tilde{z}' - \tilde{z}_0'}}
\cdot\frac{y - \tilde{z}_0'}{\v{y - \tilde{z}_0'}}\Big)^2\Big),
\end{align*}
we have
\begin{align}
\v{y - \tilde{z}'} + \v{z' - z_0'} \geq \v{y - \tilde{z}_0'} 
+ \v{z' - z_0'}\Big(1- \frac{\tilde{z}' - \tilde{z}_0'}{\v{\tilde{z}' - \tilde{z}_0'}}
\cdot\frac{y - \tilde{z}_0'}{\v{y - \tilde{z}_0'}}\Big).
\label{yogen teiri}
\end{align}
Since $z_0' \in \overline{{\mathcal U}_{1}(x)}$, 
from (\ref{bound of x,y}) and (\ref{estimate for x and z' in overline{{mathcal U}_{delta}(x)}}),
it follows that 
\begin{align}
\v{z_0' - y'} \leq \v{x' - z_0'}+\v{x' - y'} \leq L 
+ \frac{a_0}{\sqrt{1-a_0^2}}\v{x_3} \leq R, 
\label{estimate of v{z_0' - y'}}
\end{align}
where $R = L+\frac{a_0}{A\sqrt{1-a_0^2}} > 0$. 
From (\ref{estimate of v{z_0' - y'}}) and (\ref{bound of x,y}), it follows that
\begin{align*}
\Big\vert\frac{\tilde{z}' - \tilde{z}_0'}{\v{\tilde{z}' - \tilde{z}_0'}}
\cdot\frac{\tilde{z}_0' - y}{\v{\tilde{z}_0' - y}}\Big\vert
&= 
\frac{\v{(z' - z_0')\cdot(z_0' - y')}}{\v{z' - z_0'}\v{\tilde{z}_0' - y}}
\leq \frac{\v{z_0' - y'}}{\sqrt{\v{z_0' - y'}^2+y_3^2}}
\leq \frac{R}{\sqrt{R^2+A^2}},
\end{align*}
since $ t \mapsto t/\sqrt{t^2+A^2} $ is monotone increasing. 
Combining this with (\ref{yogen teiri}), we obtain
\begin{align*}
\v{y - \tilde{z}'} + \v{z' - z_0'} \geq \v{y - \tilde{z}_0'} 
+ c_2\v{z' - z_0'}
\end{align*}
where $c_2 = \frac{A^2}{\sqrt{R^2+A^2}(\sqrt{R^2+A^2}+R)} > 0$.
From (\ref{another form of tilde{l}_{x, y}}), 
it follows that
\begin{align*}
\tilde{l}_{x, y}(z') &= \frac{\v{x-\tilde{z}'_0}}{\sqrt{\gamma_-}}
+\frac{\v{z_0'-z'}+\v{\tilde{z}'-y}}{\sqrt{\gamma_+}} 
\\&
\geq \frac{\v{x-\tilde{z}'_0}}{\sqrt{\gamma_-}}
+ \frac{\v{y - \tilde{z}_0'} + c_2\v{z' - z_0'}}{\sqrt{\gamma_+}}
\\&
= l_{x, y}(z_0')+\frac{c_2}{\sqrt{\gamma_+}}\v{z' - z_0'}.
\end{align*}
Since $x \in \R^3_-$ and  
$z_0' \in \overline{{\mathcal U}_1(x)\setminus{\mathcal U}_{\delta_1}(x)}$, 
(\ref{estimate for the l_{x, y}(z') part}) implies
$$
l_{x, y}(z_0') \geq l(x, y) +c_1\v{z_0'- z'(x, y)}.
$$ 
Combining these estimates and taking $c_0 = \min\{c_1, \frac{c_2}{\sqrt{\gamma_+}}\}$, we obtain
\begin{align*}
\tilde{l}_{x, y}(z') &\geq l(x, y) + c_0(\v{z_0'- z'(x, y)}
+ \v{z_0'- z'})
\\&
\geq l(x, y) + c_0\v{z'- z'(x, y)}
\quad(x \in \overline{D}, y \in \overline{B}, z' \in \R^2\setminus{\mathcal U}_{1}(x)).
\end{align*}
This estimate and (\ref{estimate for the l_{x, y}(z') part}) imply 
(\ref{estimate of tilde{l}_{x, y}}), which 
completes the proof of Lemma \ref{the shortest length in the case}.
$\blacksquare$
\vskip1pc
%

\begin{figure}[h]
\unitlength 0.1in
\begin{picture}( 42.6000, 13.2500)(  6.5000,-26.4500)
%
\special{pn 8}%
\special{pa 810 1810}%
\special{pa 4910 1810}%
\special{fp}%
\special{sh 1}%
\special{pa 4910 1810}%
\special{pa 4844 1790}%
\special{pa 4858 1810}%
\special{pa 4844 1830}%
\special{pa 4910 1810}%
\special{fp}%
%
\special{pn 8}%
\special{pa 2810 1800}%
\special{pa 2810 2610}%
\special{da 0.070}%
\put(28.1000,-27.3000){\makebox(0,0){$x$}}%
\put(28.2000,-16.8000){\makebox(0,0){$x'$}}%
%
\special{pn 8}%
\special{pa 2820 2610}%
\special{pa 3610 1810}%
\special{da 0.070}%
%
\special{pn 8}%
\special{pa 2810 2610}%
\special{pa 1990 1820}%
\special{da 0.070}%
%
\special{pn 8}%
\special{pa 2810 2610}%
\special{pa 4180 1820}%
\special{dt 0.045}%
%
\special{pn 8}%
\special{pa 2810 2610}%
\special{pa 1410 1820}%
\special{dt 0.045}%
\put(51.8000,-17.7000){\makebox(0,0){$\R^2$}}%
%
\special{pn 8}%
\special{ar 2810 2610 310 310  3.6932476 4.6811491}%
%
\special{pn 8}%
\special{pa 4190 1800}%
\special{pa 4190 2610}%
\special{dt 0.045}%
%
\special{pn 8}%
\special{ar 4190 1810 300 300  1.5707963 2.6011732}%
\put(40.0000,-22.2000){\makebox(0,0){$\theta_0$}}%
\put(25.9000,-22.3000){\makebox(0,0){$\theta_0$}}%
\put(23.2000,-16.7500){\makebox(0,0){${\mathcal U}_{\delta_0}(x)$ }}%
\put(16.4000,-16.7500){\makebox(0,0){${\mathcal U}_{1}(x)$ }}%
%
\special{pn 25}%
\special{pa 2810 2600}%
\special{pa 4180 1810}%
\special{fp}%
\special{sh 1}%
\special{pa 4180 1810}%
\special{pa 4112 1826}%
\special{pa 4134 1838}%
\special{pa 4132 1862}%
\special{pa 4180 1810}%
\special{fp}%
%
\special{pn 25}%
\special{pa 4190 1810}%
\special{pa 4190 1810}%
\special{fp}%
%
\special{pn 25}%
\special{pa 4190 1810}%
\special{pa 4640 1810}%
\special{fp}%
\special{sh 1}%
\special{pa 4640 1810}%
\special{pa 4574 1790}%
\special{pa 4588 1810}%
\special{pa 4574 1830}%
\special{pa 4640 1810}%
\special{fp}%
%
\special{pn 25}%
\special{pa 4630 1800}%
\special{pa 4830 1500}%
\special{fp}%
\special{sh 1}%
\special{pa 4830 1500}%
\special{pa 4776 1544}%
\special{pa 4800 1544}%
\special{pa 4810 1568}%
\special{pa 4830 1500}%
\special{fp}%
%
\special{pn 25}%
\special{pa 2800 2590}%
\special{pa 3370 1810}%
\special{fp}%
\special{sh 1}%
\special{pa 3370 1810}%
\special{pa 3316 1852}%
\special{pa 3340 1854}%
\special{pa 3348 1876}%
\special{pa 3370 1810}%
\special{fp}%
%
\special{pn 25}%
\special{pa 3390 1820}%
\special{pa 3390 1820}%
\special{fp}%
%
\special{pn 25}%
\special{pa 3370 1820}%
\special{pa 3790 1560}%
\special{fp}%
\special{sh 1}%
\special{pa 3790 1560}%
\special{pa 3724 1578}%
\special{pa 3746 1588}%
\special{pa 3744 1612}%
\special{pa 3790 1560}%
\special{fp}%
%
\special{pn 8}%
\special{pa 4630 1820}%
\special{pa 4630 2590}%
\special{dt 0.045}%
%
\special{pn 8}%
\special{pa 2810 2610}%
\special{pa 4630 1820}%
\special{da 0.070}%
%
\special{pn 8}%
\special{ar 4630 1810 296 296  1.6373645 2.7091849}%
\put(43.9000,-21.9000){\makebox(0,0){$\theta$}}%
\put(41.5000,-16.2000){\makebox(0,0)[lb]{$\theta>\theta_0$}}%
\put(32.0000,-22.4000){\makebox(0,0){$\theta_0$}}%
%
\special{pn 8}%
\special{ar 2830 2620 450 450  4.6715953 5.7198690}%
\put(26.2000,-14.9000){\makebox(0,0)[lb]{refracted wave}}%
\put(63.0000,-20.9000){\makebox(0,0)[rb]{total reflected wave}}%
%
\special{pn 8}%
\special{pa 3240 1510}%
\special{pa 3550 1610}%
\special{dt 0.045}%
\special{sh 1}%
\special{pa 3550 1610}%
\special{pa 3494 1570}%
\special{pa 3500 1594}%
\special{pa 3480 1610}%
\special{pa 3550 1610}%
\special{fp}%
%
\special{pn 8}%
\special{pa 4900 1960}%
\special{pa 4680 1870}%
\special{dt 0.045}%
\special{sh 1}%
\special{pa 4680 1870}%
\special{pa 4734 1914}%
\special{pa 4730 1890}%
\special{pa 4750 1878}%
\special{pa 4680 1870}%
\special{fp}%
\end{picture}%
\caption{Propagation from the lower half-space}
\label{saitannkeironozu}
\end{figure}

\par\noindent
Proof of Proposition \ref{Asymptotics of the refracted part of the gradient of the FS}. 
For $0 < \delta_0 < 1$ given in Lemma \ref{the shortest length in the case}, 
take $\delta_1$ with $\delta_0 < \delta_1 < 1$ and put 
$
\varepsilon_0 = \frac{A}{6}\big(a_0\delta_1/\sqrt{1-a_0^2\delta_1^2}
-a_0\delta_0/\sqrt{1-a_0^2\delta_0^2}\big) > 0$.
First, we show 
\begin{align}
\{z' \in \R^2 \,\vert\, \v{z' - z'(x, y)} \leq 3\varepsilon_0 \,\} 
\subset {\mathcal U}_{\delta_1}(x)\qquad
(x \in \overline{D}, y \in \overline{B}).
\label{preparation of cutoff near z'(x, y)}
\end{align}
Choose any $x \in \overline{D}$ and $y \in \overline{B}$. 
Since (\ref{z'(x, y) is apart from the sets giving critical angle}) is equivalent to 
$ \tilde{z}'(x, y) \in \overline{{\mathcal U}_{\delta_0}(x)} $, noting 
(\ref{estimate for x and z' in overline{{mathcal U}_{delta}(x)}}) we have 
\begin{equation}
\v{x' - z'(x, y)} \leq \frac{a_0\delta_0}{\sqrt{1-a_0^2\delta_0^2}}\v{x_3}.
\label{location of z'(x, y)}
\end{equation}
If $z'$ satisfies $\v{z' - z'(x, y)} \leq 3\varepsilon_0 $, from (\ref{bound of x,y}) 
it follows that
$$
\v{z' - z'(x, y)} \leq 3\varepsilon_0 A^{-1}\v{x_3}
\leq \frac{\v{x_3}}{2}
\Big(\frac{a_0\delta_1}{\sqrt{1-a_0^2\delta_1^2}} 
- \frac{a_0\delta_0}{\sqrt{1-a_0^2\delta_0}}\Big). 
$$
This and (\ref{location of z'(x, y)}) imply
\begin{align*}
\v{x' - z'} &\leq \v{z' - z'(x, y)} + \v{x' -  z'(x, y)} 
\\&
\leq 
\frac{1}{2}\Big(\frac{a_0\delta_0}{\sqrt{1-a_0^2\delta_0^2}}
+ \frac{a_0\delta_1}{\sqrt{1-a_0^2\delta_1^2}}\Big)\v{x_3}
\\&
< \frac{a_0\delta_1}{\sqrt{1-a_0^2\delta_1^2}}\v{x_3},
\end{align*}
which yields (\ref{preparation of cutoff near z'(x, y)}) by noting 
(\ref{estimate for x and z' in overline{{mathcal U}_{delta}(x)}}).

\par

For $\delta_1$ given in (\ref{preparation of cutoff near z'(x, y)}), 
we divide the integral in 
(\ref{expression of the FS for x_3 < 0}) into three parts: 
\begin{align}
\nabla_x^k\Phi_\tau(x, y) = \frac{\tau}{4\pi\gamma_+}\big(I_{k, \tau}(x, y) 
+ J_{k, \tau}^{(1)}(x, y) + J_{k, \tau}^{(2)}(x, y)\big)
\qquad(k = 0, 1),
\label{three parts}
\end{align}
where 
$$
I_{k, \tau}(x, y) = \int_{{\mathcal U}_{\delta_1}(x)}
\nabla_x^kE^{\gamma_-}_{\tau}(x, z')
\frac{e^{-\tau\v{\tilde{z}'- y}/\sqrt{\gamma_+}}}{\v{\tilde{z}'- y}}dz',
$$
and $J_{k, \tau}^{(1)}(x, y)$ and $J_{k, \tau}^{(2)}(x, y)$ are the integrals defined by replacing
the integrated region ${\mathcal U}_{\delta_1}(x)$ in $I_{k, \tau}(x, y)$ with 
$\overline{{\mathcal U}_{1}(x)\setminus{\mathcal U}_{\delta_1}(x)}$ and 
$\R^2\setminus{\mathcal U}_{1}(x)$ respectively. Taking a cutoff function
$\phi \in C^\infty_0(\R^2)$ with $0 \leq \phi \leq 1$, $\phi(z') = 1$ 
for $\v{z'} \leq \varepsilon_0$
and $\phi(z') = 0$ for $\v{z'} \geq 2\varepsilon_0$, we define 
\begin{align*}
I_{k, \tau}^{(0)}(x, y) &= \int_{{\mathcal U}_{\delta_1}(x)}\phi(z' - z'(x, y))
\nabla_x^kE^{\gamma_-}_{\tau}(x, z')
\frac{e^{-\tau\v{\tilde{z}'- y}/\sqrt{\gamma_+}}}{\v{\tilde{z}'- y}}dz', 
\quad
\\
I_{k, \tau}^{(-\infty)}(x, y) &= \int_{{\mathcal U}_{\delta_1}(x)}
(1-\phi(z' - z'(x, y)))
\nabla_x^kE^{\gamma_-}_{\tau}(x, z')
\frac{e^{-\tau\v{\tilde{z}'- y}/\sqrt{\gamma_+}}}{\v{\tilde{z}'- y}}dz'.
\end{align*}

\par

Note that Lemma \ref{Asymptotics for no total reflection waves} and 
Proposition \ref{Asymptotics for near the total reflection angle case 1} imply 
that there exists a constant $C > 0$ such that
\begin{align}
\v{\nabla_x^kE^{\gamma_-}_{\tau}(x, z')} \leq C\tau^{k}
e^{-{\tau\v{x - \tilde{z}'}/\sqrt{\gamma_-}}}
\qquad(x \in \overline{D}, 
z' \in \overline{{\mathcal U}_{1}(x)}, k= 0, 1).
\label{simple estimate of nabla_x^kE^{gamma_-}_{tau}(x, z')}
\end{align}
Hence, for $y \in \overline{B}$ and $k = 0, 1$, from (\ref{bound of x,y}) and
(\ref{estimate of tilde{l}_{x, y}}), it follows that 
\begin{align*}
\v{J^{(1)}_{k, \tau}(x, y)} &\leq 
\int_{\overline{{\mathcal U}_{1}(x)\setminus{\mathcal U}_{\delta_1}(x)}}
\v{\nabla_xE^{\gamma_-}_{\tau}(x, z')}
\frac{e^{-\tau\v{\tilde{z}'- y}/\sqrt{\gamma_+}}}{\v{\tilde{z}'- y}}dz'
\nonumber
\\&
\leq C\tau^{k}\int_{\overline{{\mathcal U}_{1}(x)\setminus{\mathcal U}_{\delta_1}(x)}}
e^{-{\tau\v{x - \tilde{z}'}/\sqrt{\gamma_-}}}e^{-\tau\v{\tilde{z}'- y}/\sqrt{\gamma_+}}dz'
\nonumber
\\&
= C\tau^{k}\int_{\overline{{\mathcal U}_{1}(x)\setminus{\mathcal U}_{\delta_1}(x)}}
e^{-{\tau}\tilde{l}_{x, y}(z')}dz'
\nonumber
\\&
\leq C\tau^{k}e^{-{\tau}l(x, y)}
\int_{\overline{{\mathcal U}_{1}(x)\setminus{\mathcal U}_{\delta_1}(x)}}
e^{-c_0{\tau}\v{z' - z'(x, y)}}dz'.
\end{align*}
\par
For $J^{(2)}_{k, \tau}(x, y)$, noting 
Proposition \ref{Asymptotics for near the total reflection angle case 2}, 
(\ref{why tilde{l}_{x, y}(z') is defined as it is}) and 
(\ref{estimate of tilde{l}_{x, y}}), we have 
\begin{align*}
\v{J^{(2)}_{k, \tau}(x, y)} &\leq 
\int_{\R^2\setminus{\mathcal U}_{1}(x)}
\v{\nabla_xE^{\gamma_-}_{\tau}(x, z')}
\frac{e^{-\tau\v{\tilde{z}'- y}/\sqrt{\gamma_+}}}{\v{\tilde{z}'- y}}dz'
\nonumber
\\&
\leq C\tau^{k}\int_{\R^2\setminus{\mathcal U}_{1}(x)}
e^{-{\tau}T(\theta_0)}e^{-\tau\v{\tilde{z}'- y}/\sqrt{\gamma_+}}dz'
\nonumber
\\&
= C\tau^{k}\int_{\R^2\setminus{\mathcal U}_{1}(x)}
e^{-{\tau}\tilde{l}_{x, y}(z')}dz'
\nonumber
\\&
\leq C\tau^{k}e^{-{\tau}l(x, y)}
\int_{\R^2\setminus{\mathcal U}_{1}(x)}
e^{-c_0{\tau}\v{z' - z'(x, y)}}dz'.
\end{align*}
Thus, we obtain 
\begin{align}
\v{J^{(1)}_{k, \tau}(x, y)} + \v{J^{(2)}_{k, \tau}(x, y)} 
\leq C\tau^{k}e^{-{\tau}l(x, y)}
\int_{\R^2\setminus{\mathcal U}_{\delta_1}(x)}
e^{-c_0{\tau}\v{z' - z'(x, y)}}dz'.
\label{estimate of I_2+I_3 no.1}
\end{align}
Since (\ref{preparation of cutoff near z'(x, y)}) yields $\v{z' - z'(x, y)} \geq 3\varepsilon_0$ 
for any $x \in \overline{D}$, $y \in \overline{B}$ and 
$z' \in \R^2\setminus{\mathcal U}_{\delta_1}(x)$, it follows that
\begin{align*}
\int_{\R^2\setminus{\mathcal U}_{\delta_1}(x)}e^{-c_0{\tau}\v{z' - z'(x, y)}}dz' &\leq 
e^{-c_0\varepsilon_0\tau}\int_{\R^2}e^{-c_0{\tau}\v{z' - z'(x, y)}/2}dz'
\\&
\leq \frac{8\pi}{c_0^2{\tau}^2}e^{-c_0\varepsilon_0\tau}
\quad(x \in \overline{D}, y \in \overline{B}).
\end{align*}
This and (\ref{estimate of I_2+I_3 no.1}) imply that there exists a constant $C_2 > 0$
such that 
\begin{align}
\v{J^{(1)}_{k, \tau}(x, y)} + \v{J^{(2)}_{k, \tau}(x, y)} 
&\leq C_2\tau^{k-2}e^{-c_0\varepsilon_0{\tau}}e^{-{\tau}l(x, y)}
\quad(x \in \overline{D}, y \in \overline{B}, k = 0, 1).
\label{estimate of I_2 no.2}
\end{align}

\par

Since the set $\{ (x, y, z') \in \overline{D}\times\overline{B}\times\R^2 \,\vert\, 
z' \in \overline{{\mathcal U}_{\delta_1}(x)}, \v{z' - z'(x, y)} \geq \varepsilon_0 \,\}$ 
is compact, it follows that there exists a constant $c_1 > 0$ such that
$$
\frac{l_{x, y}(z') - l(x, y)}{\v{z' - z'(x, y)}} \geq 2c_1
\quad(x \in \overline{D}, y \in \overline{B}, z' \in \overline{{\mathcal U}_{\delta_1}(x)}, 
\phi(z' - z'(x, y)) \neq 1).
$$
This and (\ref{simple estimate of nabla_x^kE^{gamma_-}_{tau}(x, z')}) imply
\begin{align*}
\v{I_{k, \tau}^{(-\infty)}(x, y)} \leq C\tau^{k}e^{-{\tau}l(x, y)}
\int_{\{ z' \in \,\,{\mathcal U}_{\delta_1}(x) \,\vert\, \v{z' - z'(x, y)} \geq \varepsilon_0 \}}
e^{-2c_1{\tau}\v{z' - z'(x, y)}}dz'.
\end{align*}
Hence, similarly to getting (\ref{estimate of I_2 no.2}), we obtain
\begin{align}
\v{I_{k, \tau}^{(-\infty)}(x, y)} \leq C\tau^{k-2}e^{-c_1\varepsilon_0 \tau}e^{-{\tau}l(x, y)}
\quad(x \in \overline{D}, y \in \overline{B}, k = 0, 1).
\label{estimate of I_1^{-infty}}
\end{align}

\par

For $I_{0, \tau}^{(0)}(x, y)$, from Lemma \ref{Asymptotics for no total reflection waves}, 
it follows that
\begin{align}
I_{0, \tau}^{(0)}(x, y)
&= \frac{1}{4\pi\gamma_-}\sum_{j = 0}^{N-1}\tau^{-j}
\int_{{\mathcal U}_{\delta_1}(x)}e^{-{\tau}l_{x, y}(z')}f_{j}(z'; x, y)\phi(z'-z'(x, y))dz'
\label{the kernel representation of the FS}
\\
&
\,\,\,
+ \frac{1}{4\pi\gamma_-}
\int_{{\mathcal U}_{\delta_1}(x)}\frac{e^{-{\tau}l_{x, y}(z')}}{\v{x - \tilde{z}'}\v{\tilde{z}'-y}}
\tilde{E}_{N}(x, z'; \tau)\phi(z'-z'(x, y))dz', 
\nonumber
\end{align}
where 
\begin{align*}
f_{j}(z'; x, y) &= 
\frac{\gamma_-^{j/2}}{\v{x - \tilde{z}'}^{j+1}\v{\tilde{z}'-y}}
E_{j}(x - \tilde{z}')
\quad
(j = 0, 1, \ldots).
\end{align*}
From Lemmas \ref{Asymptotics for no total reflection waves} and 
\ref{the shortest length in the case}, the integral $I_N$ containing the remainder term 
$\tilde{E}_{N}(x, z'; \tau)$ in (\ref{the kernel representation of the FS}) is estimated by 
$$
\v{I_N} \leq C_{N, \delta_1}\frac{e^{-{\tau}l(x, y)}}{\tau^N}
\int_{\R^2}\frac{dz'}{\v{x - \tilde{z}'}^{N+1}\v{\tilde{z}'-y}}
\leq C_{N, \delta_1}\frac{e^{-{\tau}l(x, y)}}{\tau^N}
\quad(x \in \overline{D}, y \in \overline{B}, \tau \geq 1).
$$
From (\ref{preparation of cutoff near z'(x, y)}) and Lemma \ref{the shortest length in the case},
we can handle the integrals containing $f_j$ in (\ref{the kernel representation of the FS})  
as in the proof of Proposition 1 in \cite{transmission No1}. Hence, 
$I_{0, \tau}^{(0)}(x, y)$ has the same asymptotic expansion as given in 
Proposition \ref{Asymptotics of the refracted part of the gradient of the FS}. 
Similarly, we can treat $I_{1, \tau}^{(0)}(x, y)$. 
Combining these facts with (\ref{three parts}), (\ref{estimate of I_2 no.2}) 
and (\ref{estimate of I_1^{-infty}}),
we obtain Proposition \ref{Asymptotics of the refracted part of the gradient of the FS}. 
\hfill$\blacksquare$

\vskip2pc

\centerline{{\bf Acknowledgements}}

MI was partially supported by 
JSPS KAKENHI Grant Number JP17K05331. 
MK was partially supported by 
JSPS KAKENHI Grant Number JP16K05232.
This work was also partly supported by 
the Research Institute for Mathematical Sciences, 
a Joint Usage/Research Center located in Kyoto University.

%
%
%

%
%
%
\end{document}